\documentclass[reqno,12pt]{amsart}
\usepackage{a4wide}
\usepackage{mathrsfs}

\numberwithin{equation}{section}

\newtheorem{theo}{Theorem}
\newtheorem{conj}{Conjecture}
\newtheorem{coro}{Corollary}
\newtheorem{prop}{Proposition}

\theoremstyle{remark}
\newtheorem*{remark}{Remark}
\newtheorem*{remarks}{Remarks}

\def\al{\alpha}
\def\be{\beta}
\def\ga{\gamma}

\def\de{\delta}

\def\Ga{\Gamma}
\def\De{\Delta}

\def\({\left(}
\def\){\right)}
\def\[{\left[}
\def\]{\right]}

\def\dd{\textup{d}}

\date{}

\begin{document}

\title[]{How can we escape  Thomae's relations?}
\author[]{Christian Krattenthaler$^\dagger$ and Tanguy Rivoal}
\address{Institut Girard Desargues,
Universit{\'e} Claude Ber\-nard Lyon-I,
21, avenue Claude Ber\-nard,
F-69622 Villeurbanne Cedex, France}
\email{kratt@euler.univ-lyon1.fr}
\address{Institut Fourier, 
CNRS UMR 5582 / Universit\'e Grenoble 1, 
100 rue des Maths, BP 74, 
38402 Saint-Martin d'H\`eres cedex, 
France }
\email{rivoal@ujf-grenoble.fr}
\thanks{$^\dagger$ Research partially supported by EC's IHRP Programme,
grant HPRN-CT-2001-00272, ``Algebraic Combinatorics in Europe", and by
the ``Algebraic Combinatorics" Programme during Spring 2005 
of the Institut Mittag--Leffler of the Royal Swedish Academy of Sciences.}

\begin{abstract}
In 1879, Thomae discussed the relations between two generic 
hypergeometric $_3F_2$-series with argument $1$. 
It is well-known since then
that there are 120 such relations (including the trivial ones
which come from permutations of the parameters of the hypergeometric
series). More recently, Rhin and Viola
asked the following question (in a different, but equivalent language
of integrals): 
If there exists a linear dependence relation over $\mathbf{Q}$ 
between two convergent ${} _{3}  F _{2}$-series with argument $1$,
with integral parameters, and whose values are irrational numbers, 
is this relation a specialisation of one of the 120 Thomae
relations?
A few years later, Sato answered this question in the negative,
by giving six examples of relations which cannot be explained by
Thomae's relations. 
We show that Sato's counter-examples can be naturally embedded into
two families of infinitely many $_3F_2$-relations, both parametrised by
three independent parameters. 
Moreover, we find two more infinite families of the
same nature. The families, which do not seem to have been
recorded before, come from 
certain $_3F_2$-transformation formulae and contiguous relations. 
We also explain in detail the relationship between the integrals of
Rhin and Viola and $_3F_2$-series.
\end{abstract}

\maketitle

\section{Prelude: introduction and summary of the results}

In this article, we are interested in two families
of two-term relationships between  
hypergeometric ${}_{3} F_{2}$-series
with argument $1$, and the possible links between them. 
The first family consists of 120 relations found by 
Thomae~\cite{thomae}, which can be interpreted as 
the action of the symmetric group $\mathfrak{S}_5$ on five 
parameters related to the parameters of a generic ${}_{3} F_{2}$-series. 
This action has been discovered and rediscovered many times. 
We shall start our article by describing two of its seemingly different 
incarnations: one involving series (Thomae, Whipple, Hardy and others: see 
Section~\ref{sec:thomaerelations}, in particular
Theorem~\ref{prop:hardydixonandcie}) and the other 
involving integrals (Dixon, Rhin--Viola: see 
Section~\ref{sec:rhinviolagroup}, in particular
Theorem~\ref{prop:dixon}), while in Section~\ref{sec:reformulation} we
explain their equivalence. 

Our main aim is to find a hypergeometric explanation of a second family of 
six ``exotic'' integral relations recently discovered by
Sato~\cite{sato} 
(see Theorem~\ref{theo:sato}).
The latter provide counter-examples to a conjecture of Rhin and
Viola~\cite{rv} (%
see Conjecture~\ref{conj:1} in Section~\ref{sec:rhinviolagroup}),
which essentially 
predicted 
the universality of Thomae's relations in the case of integral
parameters. 
As we shall show,
this explanation is given by the following two identities
in Theorems~\ref{theo:??} and \ref{thm:gamma2}, respectively, which 
seemingly have not been stated explicitly before. 

The first one (with proof in Section~\ref{sec:escapethomae})
covers five of Sato's six original relations, and 
we will obtain from it
infinitely many explicit counter-examples to 
the conjecture by Rhin and Viola
(see Theorem~\ref{theo:2}  in Section~\ref{sec:rhinviolagroup}).

\begin{theo}\label{theo:??} Let $\al, \be, \ga$ be complex numbers
such that 
$2\al+\be+1$ and $2\be+\al+1$ are not non-positive integers, and 
such that $\Re(2\al+2\be-\ga)>0$. 
Then
\begin{equation}\label{eq:exotique2}
{} _{3} F _{2} \!\left [ \begin{matrix} {\al+1, \;\,\be+1, \;\,\gamma}
\\ {2\al+\be+1, \,2\be+\al+1}\end{matrix} 
;1\right ] = \frac{ 2(\al+\be)}{2(\al+\be)-\gamma}\,
{} _{3} F _{2} \!\left [ \begin{matrix} {\al, \quad \be, \quad\gamma}
\\ {2\al+\be+1, \,2\be+\al+1}\end{matrix} ;1\right ].
\end{equation}
\end{theo}

The second one (with proof in Section~\ref{sec:contig}) 
covers the remaining counter-example of Sato.
It implies another set of infinitely many counter-examples to
the conjecture by Rhin and Viola
(see Theorem~\ref{theo:2a}  in Section~\ref{sec:rhinviolagroup}).

\begin{theo} \label{thm:gamma2}
For any complex numbers $\al,\be,\ga$ such that  
$\Re\big(2-\be-  \frac{\al (  \al - \ga+1 ) }{\be-1}\big)>0$, and such that
$\al+1$ and $ \ga+     \frac{\al (  \al - \ga+1 ) }{\be-1}$ are
not non-positive integers, we have the identity 
\begin{multline} \label{eq:b12} 
    {} _{3} F _{2} \!\left [ \begin{matrix} { \al, \be, \ga}\\ { \al+1,\ga+ 
     \frac{\al ( \al - \ga +1) }{\be-1}+1 }\end{matrix} ;
     {\displaystyle 1}\right ]\\=
\frac { (  \al - \be +2)  ( \al + \al^2 - \ga - \al \ga + \be \ga ) } 
 {( \al+1 )  ( 2 \al + \al^2 - \al \be - \ga - \al \ga + \be \ga )  }\,
    {} _{3} F _{2} \!\left [ \begin{matrix} { \al+1, \be-1, \ga}\\ { \al+2,
\ga+     \frac{\al (  \al - \ga+1 ) }{\be-1} }\end{matrix} ;
     {\displaystyle 1}\right ].
\end{multline}
\end{theo}

Clearly, since Sato's counter-examples are special cases of 
\eqref{eq:exotique2} and \eqref{eq:b12},
but are not consequences of Thomae's relations (see Section~\ref{sec:eff}),
the two identities provide an answer to the question 
in the title. (Let us point out that Theorems~\ref{theo:??} and
\ref{thm:gamma2} are ``independent" of each other, that is, neither is it
possible to derive Theorem~\ref{thm:gamma2} from a combination of 
Theorem~\ref{theo:??} with Thomae's relations, nor is this possible in
the other direction.)

Of course, there may exist many more
ways of escaping Thomae's relations.
For example,
a rather simple-minded one consists in examining for which
integral values of the parameters a ${} _{3} F _{2}$-series with
argument $1$ can be a rational number. In fact, a complete
characterisation for the latter problem is available, see
Theorem~\ref{theo:3} in Section~\ref{sec:simple}.
Leaving this simple possibility aside,
in our proofs of identities \eqref{eq:exotique2} and \eqref{eq:b12}
in Sections~\ref{sec:escapethomae} and \ref{sec:contig}, 
we make use of two fundamentally different ways to escape
Thomae's relations:

\medskip
(1) One applies a
transformation formula transforming a $_3F_2$-series with
argument $1$ into a hypergeometric series with a larger number of
parameters (in our case, this is the transformation formula
\eqref{eq:T3240} transforming a $_3F_2$-series into 
a very-well-poised $_7F_6$-series) in order to ``exit" the
``$_3F_2$-domain," and then one ``re-enters" the ``$_3F_2$-domain" in
a different way (in our case, we use the same transformation formula
in the other direction, but after a permutation of the parameters of
the $_7F_6$-series has been carried out before). 

\medskip
(2) One starts with a
$_3F_2$-series in which one lower parameter exceeds one upper
parameter by a positive integer. Subsequently, one applies contiguous
relations to obtain a sum of several series, in which for all but one
the use of the contiguous relations has made these two parameters
equal, and thus these $_3F_2$-series with argument $1$ reduce to a
$_2F_1$-series (with argument $1$), which can then be summed by means
of the Gau{\ss} summation formula \eqref{eq:Gausz}. The various
results of these evaluations are then combined into one expression,
thereby generating a (possible huge) polynomial term, which is 
then equated to zero. In order to make this work, 
this polynomial 
must have integral solutions. (See the Remark after the proof
of Proposition~\ref{prop:gamma2} in
Section~\ref{sec:contig} for more precise explanations, and, in
particular, for an explanation of the term ``contiguous relation").

\medskip
Whereas we failed to find results other than Theorem~\ref{theo:??} by
using recipe~(1), we show in Section~\ref{sec:contig2} that recipe~(2)
can be used in many more ways than the one yielding
Theorem~\ref{thm:gamma2} (see Theorems~\ref{thm:beta} and
\ref{thm:gamma1}), thus producing many more counter-examples to
the conjecture by Rhin and Viola. 
In fact, there are certainly many more relations that can be found in
that way. We report on a curious phenomenon
in that context at the end of the ``round-up"
Section~\ref{sec:post}, where we indicate the ideas that we used
to find the hypergeometric results in Theorems~\ref{theo:??}, \ref{thm:gamma2},
\ref{thm:beta} and \ref{thm:gamma1}.

So, in summary, as disappointing as this may be, our results show
that the conjecture of Rhin and Viola was over-optimistic. The
counter-examples by Sato are not just rare exceptions, they even embed
in infinite families of counter-examples, 
and there are others beyond that. In view of
this, and since the data that we produced do not give much guidance,
we better refrain from coming up with a modified conjecture towards a
generating set of transformations for the relations 
between $_3 F_2$-series that would correct the conjecture by Rhin and Viola.
Nevertheless, finding one appears to be an interesting, and
challenging, problem.

\bigskip
\noindent {\em Acknowledgment.}
This work began during 
the second author's visit to Nihon University, Tokyo, in October 2004.
He would like to thank Noriko Hirata--Kohno for her invitation 
and Takayuki Oda who kindly gave him a copy of Sato's Master Thesis.
The first author would like to thank Anders Bj\"orner and Richard
Stanley, and the Institut Mittag--Leffler, for inviting him
to work in a relaxed and inspiring atmosphere during the
``Algebraic Combinatorics" programme in Spring 2005 at the Institut, 
during which this article was 
completed.

\section{Thomae's relations}\label{sec:thomaerelations}

Hypergeometric series are defined by
\begin{eqnarray}
{} _{q+1} F _{q} \!\left [ \begin{matrix} {\alpha_0, \alpha_1, \ldots,\alpha_{q}}\\ { 
\beta_1, \ldots, \beta_q}\end{matrix} ; {\displaystyle z}\right ] =
\sum_{k=0}^{\infty} \frac{(\al_0)_k\,(\al_1)_k\cdots (\al_q)_k}
{k!\,(\be_1)_k\cdots (\be_{q})_k} \,z^k, 
\label{eq:hyper}
\end{eqnarray} where $(\al)_0=1$ and 
$(\al)_n=\al(\al+1)\cdots (\al+n-1)$ for $n\ge 1$.
The series converges provided that 
the argument $z$ is a complex number with $\vert z \vert <1$, 
$\alpha_j\in\mathbf{C}$ and
$\beta_j\in\mathbf{C}\setminus\mathbf{Z}_{\le 0}$;
it also converges for $z= 1$ if in addition
$\Re(\be_1+\dots+\be_q)>\Re(\al_0+\dots+\al_q)$. 
Any ``permutation'' in $\mathfrak{S}_{q+1}\times\mathfrak{S}_q$ 
acting on the upper parameters 
$\al_i$, $i=0,1,\dots,q$, and the lower parameters
$\be_i$, $i=1,2,\dots,q$, 
on the left-hand side of~\eqref{eq:hyper} does not affect the
value of the right-hand side: we 
use the term ``trivial symmetries'' to indicate this fact. 

As mentioned in the introduction, 
the symmetric group $\mathfrak{S}_5$ acts classically 
on the hypergeometric series 
${}_{3} F_{2}$-series 
with argument $1$, which leads to exactly
120 formal relations between them. This group action is 
obtained using the following fundamental identity, 
due to Thomae \cite[Eq.~(12)]{thomae} (given as (3.2.2) in~\cite{Bailey}), 
which is valid under certain conditions on the parameters 
to ensure convergence of the involved series, 
\begin{align}
{} _{3} F _{2} \!\left [ \begin{matrix} {a, b, c}
\\ {d, e}\end{matrix};1
\right ] &=\frac {\Gamma(e)\,\Gamma(d+e-a-b-c)} 
{\Gamma(e-a)\,\Gamma( d+e-b-c)}
{} _{3} F _{2} \!\left [ \begin{matrix} {a, d-b, d-c}
\\ {d, d+e-b-c}\end{matrix};1
\right ].
\label{eq:thomae1}
\end{align}
The iterative application of~\eqref{eq:thomae1}, together with the
trivial symmetries, yields 120 relations, of which only 10 are inequivalent
modulo the trivial symmetries. These were given by
Thomae~\cite[Art.~4]{thomae} and put in a more suitable form by
Whipple~\cite{whipple}. 
It is apparently Hardy~\cite[p.~499]{hardy} who first gave a group
theoretic interpretation: we state his observation in the striking form given 
in~\cite{kratrao, jr}. 
\begin{theo}[\sc Hardy]
\label{prop:hardydixonandcie}
Let $s=s(x_1,x_2,x_3,x_4,x_5)=x_1+x_2+x_3-x_4-x_5$. 
The function 
\begin{equation}
\label{eq:hardy?}
\frac{1}{\Ga(s)\,\Ga(2x_4)\,\Ga(2x_5)}\,
{} _{3} F _{2} \!\left [ \begin{matrix} {2x_1-s, 2x_2-s, 2x_3-s}
\\ {2x_4,\quad 2x_5}\end{matrix};1
\right ]
\end{equation}
is a symmetric function of the five variables 
$x_1, x_2, x_3, x_4, x_5$. 
\end{theo}
Care is needed using this theorem, since $s$ is not a symmetric
function of $x_1, x_2, x_3, x_4, x_5$ and some of the ${} _{3} F
_{2}$-series might not be convergent. This result is surprising since
one could not expect {\em a priori} a much bigger invariance group than
$\mathfrak{S}_3\times\mathfrak{S}_2$, obtained by the permutations  
of $\{x_1,x_2,x_3\}$ and $\{x_4,x_5\}$, which trivially
leave~\eqref{eq:hardy?} invariant.

\section{The Rhin--Viola group for $\zeta(2)$}\label{sec:rhinviolagroup}

In 1996, Rhin and Viola introduced in~\cite{rv} the integral
\begin{equation}
\label{eq:I}
I(h,i,j,k,l)=\int_0^1 \!\!\int_0^1
\frac{x^h(1-x)^iy^k(1-y)^j}{(1-xy)^{i+j-l+1}}\,\dd x\, \dd y, 
\end{equation}
which is convergent under the assumption  
that $h,i,j,k,l$ are non-negative
integers, which will be the case throughout 
the rest of this article unless otherwise stated. Their motivation was to use the
fact that $I(h,i,j,k,l)\in\mathbf{Q} +\mathbf{Q}\zeta(2)$ to
get a good irrationality measure for 
$\zeta(2)=\sum_{n\ge 1}1/n^2=\pi^2/6$, as had been done in previous work
using similar but less general integrals (see the bibliography in~\cite{rv}).
They developed a beautiful new algebraic method for handling the
general case above and were rewarded with the best known irrationality measure
for $\pi^2$. See also~\cite{fis, zud} for related work.

{}From now on, we focus essentially 
on the hypergeometric structure underlying
their method, which is made transparent by the identity 
(see Section~\ref{sec:reformulation}, where we recall the proof for $z=1$)
\begin{equation}\label{eq:lossofsymmetry}
{} _{3} F _{2} \!\left [ \begin{matrix} {a, b, c}
\\ {d,e}\end{matrix} ; {\displaystyle z}
\right ] =  
\frac {\Gamma(d)\,\Gamma(e)}
{\Gamma(a)\,\Gamma(d-a)\,\Gamma(b)\,\Gamma(e-b)} 
\int_0^1\!\!\int_0^1
\frac{x^{a-1}(1-x)^{d-a-1}y^{b-1}(1-y)^{e-b-1}}{(1-zxy)^{c}} \,\dd x\, \dd y,
\end{equation}
which is valid provided $\Re(d)>\Re(a)>0$ and $\Re(e)>\Re(b)>0$ 
if $\vert z\vert <1$, with the further 
assumption that $\Re(d+e-a-b-c)>0$ if $z=1$. 
To simplify, we
set $B(h,i,j,k,l)=I(h,i,j,k,l)/(h!\,i!\,j!\,k!\,l!)$. The main new idea 
in~\cite{rv} was to use the action of a group on the
parameters of $h,i,j,k,l$ leaving the value of $B(h,i,j,k,l)$ invariant. To do this, 
Rhin and Viola showed that, 
under the two changes of variables $\{X=y, Y=x\}$ and $\{X=(1-x)/(1-xy),
Y=1-xy\}$, the value of $I(h,i,j,k,l)$
(and hence also that of $B(h,i,j,k,l)$) is not 
changed if the parameters are permuted by the
product of transpositions $\sigma=(h\;k)(i\;j)$ and the 5-cycle
$\tau=(h\;i\;j\;k\;l)$. The group ${\bf T}=\langle \sigma, \tau\rangle$ generated
by $\sigma$ and $\tau$ is
isomorphic to $\mathscr{D}_5$, the dihedral group of order 10: for a visual
proof, place the letters $h,i,j,k,l$, in this order, 
at the vertices of a regular pentagon.

But a more important invariance group can be obtained 
by extending 
the action of $\sigma$ and $\tau$ by linearity to the set 
$\mathscr{P}=\{h,i,j,k,l, j+k-h, k+l-i, l+h-j, h+i-k, i+j-l\}$ 
(by ``linearity'', 
we mean $\tau(h+i-k)=\tau(h)+\tau(i)-\tau(k)=i+j-l$, etc.). 
Provided  the 
five values 
$j+k-h$, $k+l-i$,  $l+h-j$, $h+i-k$, $i+j-l$
are 
non-negative (see Theorem~\ref{theo:3} in 
Section~\ref{sec:reformulation} for the arithmetic meaning of this hypothesis).
One can then use the apparent loss of the trivial symmetries in the
parameters $a,b,c$ and $d,e$ on the right-hand side of~\eqref{eq:lossofsymmetry} to prove that 
the value of $B(h,i,j,k,l)$ is invariant under the permutation on $\mathscr{P}$ defined by
$\varphi=(h\;\, i+j-l)(i\;\,l+h-j)(j+k-h\;\,k+l-i)$. Rhin and Viola managed to prove  that the group 
${\bf \Phi}=\langle \varphi, \sigma, \tau\rangle$ acting on 
$\mathscr{P}$
and leaving the value of the associated integrals invariant 
can be viewed as the 
permutation group $\mathfrak{S}_5$ acting on the set $\{h+i,\;i+j,\;j+k,\;k+l,\;l+h\}$,
and hence has cardinality 120. This remark was first made by
Dixon~\cite{dixon}, in an even  
more general form. 
\begin{theo}[\sc Dixon]\label{prop:dixon}
Assume that the complex numbers  
$h,i,j,k,l, j+k-h, k+l-i, l+h-j, h+i-k, i+j-l$ have real part $>-1$. 
Then the integral 
$B(h,i,j,k,l)$ {\em(}where $x!$ is assumed to mean 
$\Gamma(x+1)$ for complex $x${\em)} is  a symmetric function of the five parameters 
$h+i$, $i+j$, $j+k$, $k+l$, $l+h.$
\end{theo}
Finally, Rhin and Viola proposed the following conjecture.

\begin{conj}[\sc Rhin--Viola]\label{conj:1} 
Let  $h,i,j,k,l, h',i',j',k',l'$ be non-negative integers.

\hangindent15pt\hangafter1
\noindent {\em(i)} If $I(h,i,j,k,l)=I(h',i',j',k',l')$, 
then there 
exists $\rho\in{\bf T}$ 
such that  $\rho(h)=h',\rho(i)=i', \rho(j)=j', \rho(k)=k'$ and $\rho(l)=l'.$

\hangindent15pt\hangafter2
\noindent {\em(ii)} Suppose furthermore that the numbers 
\begin{align}
&j+k-h,\quad k+l-i,\quad l+h-j,\quad h+i-k,\quad i+j-l\label{eq:nonnegative1}\\
&j'+k'-h',\quad k'+l'-i',\quad l'+h'-j',\quad h'+i'-k',\quad i'+j'-l'\label{eq:nonnegative2}
\end{align}
are all non-negative. If 
$I(h,i,j,k,l)/I(h',i',j',k',l')\in\mathbf{Q}$, 
then there 
exists $\rho\in{\bf \Phi}$ such that
$\rho(h)=h',\rho(i)=i', \rho(j)=j', \rho(k)=k'$ and $\rho(l)=l'.$
\end{conj}

The truth of (i) and (ii) would have shown that their
method is optimal, but both have been shown to be false in 2001 by
Susumu Sato~\cite{sato}, who found the following counter-examples,
apparently by numerical inspection.

\begin{theo} [\sc Sato] \label{theo:sato}
Both cases of Conjecture~\ref{conj:1} are false, as shown by the following 
six counter-examples:
\begin{align}
I(1,1,1,1,1) &=5-3\zeta(2)=I(3,1,1,2,0),\label{eq:sato1}\\
I(3,1,2,2,1) &=79/4-12\zeta(2)=I(4,2,2,3,0),\label{eq:sato2}\\
I(3,1,2,1,1) &=3\zeta(2)-59/12= I(3,3,1,3,0),\label{eq:sato4}\\
I(3,2,2,2,1) &=10\zeta(2)-148/9=I(5,1,3,2,1),\label{eq:sato5}\\
I(3,0,3,1,1) &=9\zeta(2)-59/4=9I(3,3,1,2,1),\label{eq:sato6}\\
I(3,1,3,1,0) &=\zeta(2)-29/18= I(3,2,1,2,0).\label{eq:sato3}
\end{align} 
\end{theo}
(Sato mis-stated~\eqref{eq:sato6}  as $I(3,0,3,1,1)=I(3,3,1,2,1)$.
The reader should also note that \eqref{eq:sato4} and \eqref{eq:sato6}
altogether relate four different integrals rationally.)
Equation~\eqref{eq:sato1} is already a counter-example to both (i) and (ii).
The following
questions are natural, but were not considered by Sato:

\begin{itemize}
\item Are these counter-examples merely numerical accidents, or 
  do they admit a theoretical explanation?

\item Do there exist infinitely many counter-examples to the
  conjecture of Rhin and Viola?
\end{itemize}
We give a complete answer to both questions in the two theorems below
which we prove in Sections~\ref{sec:proof} and \ref{sec:proofa},
respectively.

\begin{theo}\label{theo:2} 
{\em(i)} Sato's counter-examples~\eqref{eq:sato1} up to \eqref{eq:sato6} can be
  explained by purely hypergeometric means, i.e., there exists a
  general hypergeometric identity that generates them.

\noindent {\em(ii)} For each integer $\al\ge 1$, the equation
\begin{equation}\label{eq:infinicontreexemples}
I(2\al-1,2\al-1,\al, 2\al-1,\al)=I(2\al+1,2\al-1,\al, 2\al,\al-1)
\end{equation}
provides a counter-example to the
cases {\em(i)} and {\em(ii)} of Conjecture~{\em\ref{conj:1}}.
\end{theo}

\begin{theo}\label{theo:2a} 
{\em(i)} Sato's counter-example~\eqref{eq:sato3} can be
  explained by purely hypergeometric means, i.e., there exists a
  general hypergeometric identity that generates them.

\noindent {\em(ii)} For each integer $\al\ge 2$, the equation
\begin{equation}\label{eq:infinicontreexemplesa}
I(\al^2-1,\al-1,\al^2-\al+1,\al-1,0)=(\al-1)\,I(\al^2-1,\al,\al^2-\al-1,\al,0)
\end{equation}
provides a counter-example to 
case {\em(ii)} of Conjecture~{\em\ref{conj:1}}, and also
to case {\em(i)} if $\al=2$.
\end{theo}

\begin{remarks}
(1)
A particularly elegant instance of \eqref{eq:exotique2} is the one
where $\ga=\al+\be$:
for any complex numbers $\al$ and $\be$
which are not non-positive integers and which satisfy 
$\Re(\al+\be) >0$, we have 
\begin{equation}\label{eq:exotique}
{} _{3} F _{2} \!\left [ \begin{matrix} {\al+1, \;\,\be+1, \;\,\al+\be}
\\ {2\al+\be+1, 2\be+\al+1}\end{matrix};1 
\right ] = 2\,
{} _{3} F _{2} \!\left [ \begin{matrix} {\al, \quad \be, \quad\al+\be}
\\ {2\al+\be+1, 2\be+\al+1}\end{matrix};1 
\right ].
\end{equation}
The action of Thomae's relations on both sides of~\eqref{eq:exotique} 
independently provides  
ten variations of~\eqref{eq:exotique}, up to trivial symmetries: 
one example, given 
by~\eqref{eq:couplage} in Section~\ref{sec:proof}, 
will be used in the proof 
of Theorem~\ref{theo:2}. 
Equation~\eqref{eq:infinicontreexemples} follows 
from the case of~\eqref{eq:exotique} where $\al=\be$ is a positive integer. 
Furthermore, we shall show in Section~\ref{sec:proof} that 
$I(2\al-1,2\al-1,\al, 2\al-1,\al)$ tends to 0 as $\al$ tends to
infinity: this fact implies  
that~\eqref{eq:infinicontreexemples} provides infinitely many 
counter-examples to Conjecture~\ref{conj:1}. 

(2)
As we show in Section~\ref{sec:pattern},
the ``general'' hypergeometric identity that generates~\eqref{eq:sato1} 
up to \eqref{eq:sato6} is exactly identity~\eqref{eq:exotique2},
via the translation between integrals and hypergeometric
$_3F_2$-series given in \eqref{eq:lossofsymmetry}.
Equation \eqref{eq:infinicontreexemples} is a special case. 

(3)
Similarly, we show in Section~\ref{sec:pattern} that
the ``general'' hypergeometric identity that generates~\eqref{eq:sato3} 
is exactly identity~\eqref{eq:b12}, again via the translation
\eqref{eq:lossofsymmetry}.
Equation \eqref{eq:infinicontreexemplesa} is a special case. 
Since we show there that the integral on the left-hand side of
\eqref{eq:infinicontreexemplesa} tends to zero as $\al$ tends to
infinity, also \eqref{eq:infinicontreexemplesa} provides infinitely
many counter-examples to Conjecture~\ref{conj:1}.

(4)
It would also be interesting to look at the analogous problem 
arising from the group action on the triple integral 
$$
\int_0^1\!\!\int_0^1\!\!\int_0^1  \frac{u^h(1-u)^lv^k(1-v)^sw^j(1-w)^q}
{(1-(1-uv)w)^{q+h-r+1}}\,\dd u\,\dd v\,\dd w \in\mathbf{Q}+\mathbf{Q}\zeta(3)
$$
found by Rhin and Viola in~\cite{rv2}: do there exist 
exotic relationships between such integrals that are not described by this group action? 
Note that this action admits an interpretation in terms of 
very-well-poised ${} _{7} F _{6}$-series 
exactly in the style of Theorem~\ref{prop:hardydixonandcie}
(see \cite[Sec.~4]{zud} for the passage from the integrals to 
very-well-poised ${} _{7} F _{6}$-series, and 
\cite[Proposition~5, $q\to 1$, p.~6698]{jr} for a particularly elegant
formulation of the group structure). For very clear 
expositions of various group actions on $(q-)$ hypergeometric series, 
see~\cite{kratrao,jr} and 
the references therein.
\end{remarks}

\section{From Dixon to Thomae
}\label{sec:reformulation}

In this section, we show more precisely how Rhin and Viola's integrals
are related to hypergeometric series. 
To get a new expression for the integral $I(h,i,j,k,l)$, we
transform the integrand of~\eqref{eq:I} by using the 
binomial series expansion 
\begin{equation*}
\frac{1}{(1-xy)^{i+j-l+1}} =\sum_{n=0}^{\infty}
\frac{(i+j-l+1)_n}{n!}\, (xy)^{n},
\end{equation*}
and the beta integral evaluations
$$
\int_0^1 x^{n+h}(1-x)^{i}\dd x = \frac{(n+h)!\,i!}{(n+h+i+1)!},\quad 
\int_0^1 y^{n+k}(1-y)^{j}\dd y = \frac{(n+k)!\,j!}{(n+k+j+1)!}.
$$
We have
\begin{align}
B(h,i,j,k,l)
&=\frac{1}{h!\,k!\,l!}\,\sum_{n=0}^{\infty}
\frac{(n+h)!\,(n+k)!\,(i+j-l+1)_n}{n!\,(n+h+i+1)!\,(n+k+j+1)!}\label{eq:fin4}\\
&=\frac{1}{l!\,(h+i+1)!\,(k+j+1)!}\,
{} _{3} F _{2} \!\left [ \begin{matrix} {h+1, k+1, i+j-l+1}
\\ {h+i+2,k+j+2}\end{matrix};1\right]\label{eq:fin3},
\end{align}
since the interchange of summation and integral is justified by 
Fubini's theorem. 
The passage from~\eqref{eq:fin4} to~\eqref{eq:fin3} uses 
the trivial identity $(\al+n)!=\al!\,(\al+1)_n$.

Under this interpretation, it is not surprising that 
the group obtained by Dixon and Rhin--Viola should be a reformulation 
of Theorem~\ref{prop:hardydixonandcie}, in 
terms of integrals rather than series. 
Indeed, if we define a bijection between the tuples 
$(x_1, x_2, x_3, x_4, x_5)$ and $(h,i,j,k,l)$ by 
\begin{align*}
h+1&=\hphantom{-}x_1-x_2-x_3+x_4+x_5,\\
i+1&=-x_1+x_2+x_3+x_4-x_5,\\
j+1&=\hphantom{-}x_1-x_2+x_3-x_4+x_5,\\
k+1&=-x_1+x_2-x_3+x_4+x_5,\\
l+1&=\hphantom{-}x_1+x_2+x_3-x_4-x_5=s,
\end{align*}
then we see that 
$B(h,i,j,k,l)$, written as~\eqref{eq:fin3},  perfectly matches~\eqref{eq:hardy?} 
and the 120 possible series are all convergent if the ten integers in 
the set $\mathscr{P}$ are non-negative. 
Since 
\begin{equation}
\label{eq:xhijkl}
2x_1=l+h+2,\;
2x_2=k+l+2,\;
2x_3=i+j+2,\;
2x_4=h+i+2,\;
2x_5=j+k+2,
\end{equation}
we also see that the symmetry of~\eqref{eq:hardy?} in the variables 
$x_1, x_2, x_3, x_4, x_5$ is 
equivalent to  
the symmetry of $B(h,i,j,k,l)$ in the variables 
$h+i$, $i+j$, $j+k$, $k+l$, $l+h.$

\section{When is $I(h,i,j,k,l)$ rational?}
\label{sec:simple}

In this section we answer the question of ``simple-minded" counter-examples
to Conjecture~\ref{conj:1} that was raised in the Introduction.

\begin{theo}\label{theo:3}
Let $h,i,j,k,l$ be non-negative integers. 
Then the following assertions are equivalent:
\begin{itemize}
\item[(a)] The integers $j+k-h$, $k+l-i$, $l+h-j$, $h+i-k$, $i+j-l$ are all non-negative.
\item[(b)] The integral $I(h,i,j,k,l)$ is an irrational number.
\end{itemize}
\end{theo}
\begin{remarks}
(1) We remarked earlier (see \eqref{eq:lossofsymmetry} and
\eqref{eq:fin3}) that
$I(h,i,j,k,l)$ is essentially equal to a $_3F_2$-series.
If we translate Theorem~\ref{theo:3} into the analogous theorem for
\begin{equation} \label{eq:3F2} 
{} _{3} F _{2} \!\left [ \begin{matrix} a,b,c\\
d,e
\end{matrix};1  \right ],
\end{equation}
via the relations $a=h+1$, $b=k+1$, $c=i+j-l+1$, $d=h+i+2$, $e=k+j+2$,
we get the following necessary and sufficient condition for the
irrationality of the series \eqref{eq:3F2} for integral values of
$a,b,c,d,e$:
\begin{gather}
\label{eq:cond1}
d+e\ge a+b+c+1,\quad a\ge1,\quad b\ge1,\quad c\ge1,\\
d\ge \max\{ a, b, c\}+1,\\
e\ge \max\{ a, b, c\}+1.
\label{eq:cond2}
\end{gather}

(2) 
As a marginal consequence, Theorem~\ref{theo:3} 
proves that the analogue of the case (ii) of 
Conjecture~\ref{conj:1}, where we now suppose that 
the non-negativity condition~\eqref{eq:nonnegative1} is not true, cannot hold either. Indeed, 
if one of the integers in
\eqref{eq:nonnegative1} were negative and none in~\eqref{eq:nonnegative2} were 
negative, then the value of $I(h,i,j,k,l)/I(h',i',j',k',l')$ 
would be irrational 
and the conjecture would be empty. 
And if one of the integers in
\eqref{eq:nonnegative1} and one in~\eqref{eq:nonnegative2} were 
negative, then the conjecture would be trivially false because, although
$I(h,i,j,k,l)$ and $I(h',i',j',k',l')$ are rational, there exist many 
${\bf \Phi}$-unrelated choices for $h, h'$, etc.: one may consider 
$I(1,1,1,1,3)$ and $I(1,1,1,1,4)$ for example.
\end{remarks}

\begin{proof}[Proof of Theorem~{\em\ref{theo:3}}]
We first show the implication (b)${}\Rightarrow{}$(a). Since the parameters 
$j+k-h$, $k+l-i$, $l+h-j$, $h+i-k$, $i+j-l$ are cyclically permuted by $\tau\in{\bf T}$, if one of them 
were negative, then without loss of generality, we may 
assume that it is $i+j-l$. But $i+j-l\le -1$ implies that the 
integrand of $I(h,i,j,k,l)$ is a polynomial with integral coefficients and 
hence that $I(h,i,j,k,l)\in\mathbf{Q}.$

The reverse implication (a)${}\Rightarrow{}$(b) is a little
bit more complicated. 
Since $i ,j$ and $i+j-l$ are positive integers, 
we can write the expansion~\eqref{eq:fin4} in the equivalent form: 
\begin{equation}
\label{eq:rationalfraction}
I(h,i,j,k,l)=\frac{i!\,j!}{(i+j-l)!}\,\sum_{n=0}^{\infty}
\frac{(n+1)_{i+j-l}}{(n+h+1)_{i+1}(n+k+1)_{j+1}}.
\end{equation}
(We have used trivial identities such as 
$(n+h+i)!/(n+h)!=(n+h+1)_{i+1}$.)
We know that $I(h,i,j,k,l)\in\mathbf{Q}+\mathbf{Q}\zeta(2)$ and it will be enough to 
prove that the coefficient $p(h,i,j,k,l)$ of the irrational number $\zeta(2)$ is non-zero. A standard 
way to find an explicit expression for 
this coefficient is to expand  the summand of~\eqref{eq:rationalfraction}, which is a 
rational function of $n$, 
in partial fractions 
(see the introduction of~\cite{ri} for details in many similar 
cases and references). All computations done, one finds that
$$
p(h,i,j,k,l)=(-1)^{h+i+j+k+l}\sum_{s=\max(h,k,i+j-l)}^{\min(h+i,k+j)}
\binom{i}{s-h}\binom{j}{s-k}\binom{s}{i+j-l},
$$ 
with the convention that the value of the sum is $0$ if it is empty.
The latter is the case if and only if  
$\min(h+i,k+j)<\max(h,k,i+j-l)$.

We now  show that condition (a) ensures that the sum is non-empty and hence that 
$$
(-1)^{h+i+j+k+l}p(h,i,j,k,l)>0,
$$
because it is a sum of  binomial coefficients.
We have already used the fact that $i+j-l\ge 0$. Since the inequalities $h+i-k\ge 0$ 
and  $k+j-h\ge 0$ imply that 
$\max(h,k)\le\min(h+i,k+j)$, it only remains to show that $i+j-l\le \min(h+i,k+j)$ to finally 
prove that $\min(h+i,k+j)\ge \max(h,k,i+j-l)$. 
But $\min(h+i,k+j)-(i+j-l)=\min(h+l-j, k+l-i)\ge 0$,  which finishes 
the proof.
%
\end{proof}

\section{Effective computation of Thomae relations}
\label{sec:eff}

In this section, we show how to compute the complete set of (generically) 120 
Thomae relations (convergent or not) for 
any given ${} _{3}  F _{2}$-series with argument $1$.
We need this to transform Sato's counter-examples into more suitable 
forms. 
The most effective way to do this is by using the
parametrisation $2x_1-s,2x_2-s,2x_3-s$ and $2x_4, 2x_5$ 
of the upper and lower parameters of the $_3F_2$-series from
Theorem~\ref{prop:hardydixonandcie}. If we denote the upper parameters
by $a,b,c$ and the lower parameters by $d,e$, then 
$x_1,x_2,x_3,x_4,x_5,s$ and $a,b,c,d,e$ are related by
$$
2x_1=d+e-b-c, \, 2x_2=d+e-c-a,\, 2x_3=d+e-a-b, \, 2x_4=d,\, 2x_5=e,
$$ 
and 
$$s=x_1+x_2+x_3-x_4-x_5=2(d+e-a-b-c).$$ 
If one prefers the integral setting, then one gets the  
parametrisation of Theorem~\ref{prop:hardydixonandcie} of the integral
$I(h,i,j,k,l)$ by the formulae~\eqref{eq:xhijkl} given in 
Section~\ref{sec:reformulation}.

The following simple {\em Maple} commands  compute all possible values of 
the arrays of parameters $[2x_{\rho(1)}-s_{\rho},2x_{\rho(2)}-s_{\rho}, 
2x_{\rho(3)}-s_\rho ; 2x_{\rho(4)},2x_{\rho(1)}]$ (with $s_\rho=x_{\rho(1)}+ 
x_{\rho(2)}+x_{\rho(3)}-x_{\rho(4)}-x_{\rho(1)}$) 
over all permutations $\rho$ of $\{1,2,3,4,5\}$, with the nice feature
to output only the term-wise different arrays (viewed as 5-tuples by the
program): 
\medskip

\begin{small}
\noindent $> {\tt with(combinat)\!:}$
\newline
$> {\tt p:=(u,v,w,x,y)-\!\!>permute([u,v,w,x,y])\!:}$
\newline
$> {\tt s:=(u,v,w,x,y)-\!\!>u\!+\!v\!+\!w\!-\!x\!-\!y\!:}$
\newline
$> {\tt A:=(u,v,w,x,y)-\!\!>[2\!*\!u\!-\!s(u,v,w,x,y), 
2\!*\!v\!-\!s(u,v,w,x,y), 2\!*\!w\!-\!s(u,v,w,x,y), 2\!*\!x,2\!*\!y]\!:}$
\newline
$> {\tt T:=(u,v,w,x,y)-\!\!>\tt seq(A(op(1,op(j,p(u,v,w,x,y))),op(2,op(j,p(u,v,w,x,y))),}$
\newline
\indent ${\tt op(3,op(j,p(u,v,w,x,y))),op(4,op(j,p(u,v,w,x,y))),op(5,op(j,p(u,v,w,x,y)))),}$
\newline
\indent ${\tt j=1...nops(p(u,v,w,x,y)))\!:}$
\newline
$> {\tt F:=(a,b,c,d,e)-\!\!>T((d\!+\!e\!-\!b\!-\!c)/2,(d\!+\!e-\!c\!-\!a)/2,(d\!+\!e\!-\!a\!-\!b)/2,d/2,e/2)\!:}$
\newline
$> {\tt I:=(h,i,j,k,l)-\!\!>T((h\!+\!l\!+\!2)/2,(k\!+\!l\!+\!2)/2,(i\!+\!j\!+\!2)/2,(h\!+\!i\!+\!2)/2,
(j\!+\!k\!+\!2)/2)\!:}$
\end{small}
\medskip

\noindent 
The function {\tt T} computes all the different expressions 
for the value of the symmetric function in Theorem~\ref{prop:hardydixonandcie}, {\tt F} does the same for a 
${} _{3}  F _{2}[a,b,c;d,e]$ and $\tt I$ for (\footnote{
The letter ``{\tt I}'' denotes the complex number $i$ in {\em Maple} and one must use 
another symbol. But since \LaTeX~is not {\em Maple}, there is no problem here.}) 
$I(h,i,j,k,l)$. Only the Gamma-factors are not computed, but 
this could be easily done. For example, we obtain
\medskip

\begin{small}
\noindent $> {\tt I(1,1,1,1,1);}$
$$
{\tt [2,2,2,4,4]}
$$
$> {\tt I(3,1,1,2,0);}$
$$
{\tt [4, 3, 3, 6, 5], \;[2, 1, 3, 4, 5],\; [3, 2, 3, 6, 4], \;[1, 1, 2, 4, 4], \;[2, 4, 2, 5, 5].}
$$
\end{small}%
\noindent {\em Maple} outputs 25 other arrays for $I(3,1,1,2,0)$
 but since they correspond to the five above by the trivial symmetries, we do not list them.
We can also find the Thomae relations for both sides of
counter-example~\eqref{eq:sato3}:
\medskip

\begin{small}
\noindent $> {\tt I(3,1,3,1,0);}$
$$
\!{\tt [4, 2, 5, 6, 6], \;[1, 2, 2, 6, 3],\;  [1, 4, 4, 5, 6], 
 \;[1, 1, 1, 3, 5]}
$$
$>{\tt I(3,2,1,2,0);}$
$$
{\tt [4, 3, 4, 7, 5],\; [2, 1, 4, 5, 5],\; [1, 1, 3, 4, 5],\; [3, 3, 3, 4, 7].}
$$
\end{small}

\noindent This shows that the relations $I(1,1,1,1,1)\!=\!I(3,1,1,2,0)$ and 
$I(3,1,3,1,0)\!=\!I(3,2,1,2,0)$ are 
not consequences of Thomae relations. Similar computations 
provide a verification of the other 
counter-examples.

\section{The pattern behind Sato's counter-examples }\label{sec:pattern}

With the interpretation given in Section~\ref{sec:reformulation}, the
case (ii) of Conjecture~\ref{conj:1} can be reformulated as follows:

\medskip
{\em If there exists a linear dependence relation over $\mathbf{Q}$
between two convergent ${} _{3}  F _{2}$-series with argument $1$,
with integral parameters, and whose values are irrational numbers, 
then this relation is a specialisation of one of the 120 Thomae
relations.}
\medskip

Sato's counter-examples destroy this hope.
We can formulate his counter-examples~\eqref{eq:sato1} and~\eqref{eq:sato2}
in hypergeometric form (with simplification of the Gamma-fac\-tors) as
follows: 
\begin{equation*} 
\,{} _{3} F _{2} \!\left [ \begin{matrix} {2, 2, 2}
\\ {4, 4}\end{matrix} ;1\right] 
= \frac{3}{20} \,
{} _{3} F _{2} \!\left [ \begin{matrix} {4, 3, 3}
\\ {6, 5}\end{matrix};1 \right] 
\label{eq:3f2-1}
\quad \textup{and}\quad 
\,{} _{3} F _{2} \!\left [ \begin{matrix} {4, 3, 3}
\\ {6, 6}\end{matrix};1 \right] 
= \frac{2}{21}\,
{} _{3} F _{2} \!\left [ \begin{matrix} {5, 4, 5}
\\ {8, 7}\end{matrix} ;1\right].
\end{equation*}
Under this form, the parameters on
the left-hand sides and those on the right-hand sides seem still 
rather unrelated, and it is thus 
still unclear whether we face numerical accidents
or if there is something deeper behind. 

However, a natural thing to do here is to seek new numerical relations by
applying Thomae's transformations (using the {\em Maple} commands 
of the previous section) 
to each of the four ${} _{3} F _{2}$-series in
\eqref{eq:3f2-1}, independently. 
We find that we are trying to prove that
\begin{equation} \label{eq:counter1}
{} _{3} F _{2} \!\left [ \begin{matrix} {2, 2, 2}
\\ {4, 4}\end{matrix};1 \right] = 2\,
{} _{3} F _{2} \!\left [ \begin{matrix} {1, 1, 2}
\\ {4, 4}\end{matrix} ;1\right] \quad \textup{and}\quad 
{} _{3} F _{2} \!\left [ \begin{matrix} {3, 2, 3}
\\ {6, 5}\end{matrix};1 \right] = 2\,
{} _{3} F _{2} \!\left [ \begin{matrix} {2, 1, 3}
\\ {6, 5}\end{matrix};1 \right],
\end{equation}
where a pattern now emerges, explained by
the earlier identity~\eqref{eq:exotique}. 

The hypergeometric forms 
of the three  coun\-ter-examples~\eqref{eq:sato4},~\eqref{eq:sato5} 
and~\eqref{eq:sato6} are 
\begin{multline} \label{eq:counter2}
 {} _{3} F _{2} \!\left [ \begin{matrix} {4, 2, 3}
\\ {6, 5}\end{matrix};1 \right] = \frac{3}{35}\,
{} _{3} F _{2} \!\left [ \begin{matrix} {4, 4, 5}
\\ {8, 6}\end{matrix};1 \right],\quad 
{} _{3} F _{2} \!\left [ \begin{matrix} {4, 3, 4}
\\ {7, 6}\end{matrix};1 \right] = \frac{5}{7}\,
{} _{3} F _{2} \!\left [ \begin{matrix} {6, 3, 4}
\\ {8, 7}\end{matrix};1 \right],\\
{} _{3} F _{2} \!\left [ \begin{matrix} {4, 2, 3}
\\ {5, 6}\end{matrix};1 \right]=\frac{1}{21}\,
{} _{3} F _{2} \!\left [ \begin{matrix} {4, 3, 4}
\\ {8, 5}\end{matrix};1 \right],
\end{multline}
which 
become much more illuminating when rewritten as
\begin{multline*}
{} _{3} F _{2} \!\left [ \begin{matrix} {3, 2, 4}
\\ {6, 5}\end{matrix};1 \right] = \frac13\,
{} _{3} F _{2} \!\left [ \begin{matrix} {2, 1, 4}
\\ {6, 5}\end{matrix};1 \right], \quad 
{} _{3} F _{2} \!\left [ \begin{matrix} {3, 2, 2}
\\ {6, 5}\end{matrix};1 \right] = \frac23\,
{} _{3} F _{2} \!\left [ \begin{matrix} {2, 1, 2}
\\ {6, 5}\end{matrix};1 \right], \\
{} _{3} F _{2} \!\left [ \begin{matrix} {3, 2, 4}
\\ {6, 5}\end{matrix};1 \right] = \frac13\,
{} _{3} F _{2} \!\left [ \begin{matrix} {2, 1, 4}
\\ {6, 5}\end{matrix};1 \right],
\end{multline*}
by using Thomae's relations. (In particular, \eqref{eq:sato4} 
and~\eqref{eq:sato6} are consequences of the same identity). The
connexion with Theorem~\ref{theo:??} is now clear.

Finally,
the hypergeometric form of the 
counter-example~\eqref{eq:sato3} is
\begin{equation} \label{eq:sato3hyp} 
   {{} _{3} F _{2} \!\left [ \begin{matrix} { 4, 2, 5}\\ { 6,
      6}\end{matrix} ; {\displaystyle 1}\right ] } =\frac {5} {9}\,
{{} _{3} F _{2} \!\left [ \begin{matrix} { 4, 3, 4}\\ { 7,
      5}\end{matrix} ; {\displaystyle 1}\right ] } ,
\end{equation}
which is obviously the special case $\al=4$, $\be=3$, $\ga=4$ 
of \eqref{eq:b12}.
 
\section{Proof
of Theorem~\ref{theo:??} }\label{sec:escapethomae}

For the proof of \eqref{eq:exotique2} we need the following
transformation formula due to Verma and Jain
(see \cite[(3.5.10), $q\to 1$, reversed]{GaRaAA}, 
being implied by \cite[(4.1)]{VeJaAD}) 
between a $_3F_2$-series and a very-well-poised
$_7F_6$-series: 
\begin{multline} \label{eq:T3240}
{} _{3} F _{2} \!\left [ \begin{matrix} { b, c,d}\\ { a, a - b + c}\end{matrix} ;
   {\displaystyle 1}\right ]  =
  \frac{
      \Gamma({ \textstyle 2 a}) \,\Gamma({ \textstyle 2 a - 2 b - d}) \,
      \Gamma({ \textstyle a - b + c}) \,\Gamma({ \textstyle a - d + c}) }
      {\Gamma({ \textstyle 2 a - 2 b}) \,\Gamma({ \textstyle 2 a - d}) \,
      \Gamma({ \textstyle a + c}) \,\Gamma({ \textstyle a - b - d +
c}) }\\
\times
{} _{7} F _{6} \!\left [ \begin{matrix} { 
       a- \frac{1}{2}, \frac{a}{2}+\frac{3}{4}, b, \frac{d}{2}, 
       \frac{d}{2}+\frac{1}{2}, \frac{a}{2} - \frac{c}{2}, \frac{a}{2} -
       \frac{c}{2}+ \frac{1}{2}}\\ {  \frac{a}{2}- \frac{1}{4},
        a - b+\frac{1}{2},  a - \frac{d}{2}+\frac{1}{2}, a - \frac{d}{2},
       \frac{a}{2} + \frac{c}{2}+\frac{1}{2}, \frac{a}{2} +
       \frac{c}{2}}\end{matrix} ; {\displaystyle 1}\right ].
\end{multline}
If we apply this transformation to the $_3F_2$-series on the left-hand
side of \eqref{eq:exotique2}, then we obtain
\begin{multline} \label{eq:7F6a}
 \frac{ \Gamma({ \textstyle \al + 2 \be+1}) \,
     \Gamma({ \textstyle 4 \al + 2 \be+2}) \,
     \Gamma({ \textstyle 2 \al + 2 \be -
      \ga}) \,\Gamma({ \textstyle 2 \al +
      2 \be - \ga+2}) }{\Gamma({ \textstyle
      2 \al + 2 \be}) \,
     \Gamma({ \textstyle 2 \al + 2 \be+2}) \,
     \Gamma({ \textstyle \al + 2 \be -
      \ga+1}) \,\Gamma({ \textstyle 4 \al +
      2 \be - \ga+2}) }\\
\times
{} _{7} F _{6} \!\left [ \begin{matrix} { 2 \al +
      \be+\frac{1}{2}, \al + \frac{\be}{2}+\frac{5}{4},
      \al+1, \frac{\ga}{2}, 
      \frac{\ga}{2}+\frac{1}{2}, \al, 
      \al+\frac{1}{2}}\\ { \al +
      \frac{\be}{2}+\frac{1}{4},  \al + \be+\frac{1}{2},
      2 \al + \be -
      \frac{\ga}{2}+\frac{3}{2}, 2 \al + \be -
      \frac{\ga}{2}+1,  \al + \be+\frac{3}{2},
      \al + \be+1}\end{matrix} ; {\displaystyle 1}\right
      ] .
\end{multline}
We permute the parameters in the $_7F_6$-series to get the equivalent
expression 
\begin{multline} \label{eq:7F6b}
 \frac{ \Gamma({ \textstyle \al + 2 \be+1}) \,
     \Gamma({ \textstyle 4 \al + 2 \be+2}) \,
     \Gamma({ \textstyle 2 \al + 2 \be -
      \ga}) \,\Gamma({ \textstyle 2 \al +
      2 \be - \ga+2}) }{\Gamma({ \textstyle
      2 \al + 2 \be}) \,
     \Gamma({ \textstyle 2 \al + 2 \be+2}) \,
     \Gamma({ \textstyle \al + 2 \be -
      \ga+1}) \,\Gamma({ \textstyle 4 \al +
      2 \be - \ga+2}) }\\
\times
{} _{7} F _{6} \!\left [ \begin{matrix} {  2 \al +
      \be+\frac{1}{2},  \al + \frac{\be}{2}+\frac{5}{4},
      \al, \frac{\ga}{2}, 
      \frac{\ga}{2}+\frac{1}{2},  \al+\frac{1}{2}, 
      \al+1}\\ {  \al +
      \frac{\be}{2}+\frac{1}{4},  \al + \be+\frac{3}{2},
       2 \al + \be -
      \frac{\ga}{2}+\frac{3}{2},  2 \al + \be -
      \frac{\ga}{2}+1,  \al + \be+1,
      \al + \be+\frac{1}{2}}\end{matrix} ;
      {\displaystyle 1}\right ] .
\end{multline}
To this $_7F_6$-series, 
we apply the transformation \eqref{eq:T3240} in the backward
direction, that is we apply the transformation
\begin{multline*}
 {} _{7} F _{6} \!\left [ \begin{matrix} { a,  \frac{a}{2}+1, b, c, 
    c+\frac{1}{2}, d, d+\frac{1}{2}}\\ { \frac{a}{2}, a - b+1, a -
c+1, a - c+\frac{1}{2}, a - d+1, a - d+\frac{1}{2}}\end{matrix} ; 
{\displaystyle
    1}\right ]  \\= \frac{  
      \Gamma({ \textstyle 2 a - 2 b+1}) \,
      \Gamma({ \textstyle 2 a - 2 c+1}) \,
      \Gamma({ \textstyle 2 a - 2 d+1}) \,
      \Gamma({ \textstyle 2 a - b - 2 c - 2 d+1}) }{\Gamma({ \textstyle 
       2 a+1}) \,\Gamma({ \textstyle 2 a - 2 b - 2 c+1}) \,
      \Gamma({ \textstyle 2 a - b - 2 d+1}) \,
      \Gamma({ \textstyle 2 a - 2 c - 2 d+1}) }\\
\times
{} _{3} F _{2} \!\left [ \begin{matrix} { 2 c, b,
       a - 2 d+\frac{1}{2}}\\ { 2 a - b - 2 d+1, 
       a+\frac{1}{2}}\end{matrix} ; {\displaystyle 1}\right ].
\end{multline*}
Thus we directly arrive at the right-hand side of \eqref{eq:exotique2}.

\section{Proof of Theorem~\ref{theo:2}} \label{sec:proof}

As already mentioned in Section~\ref{sec:pattern}, the cases 
$(\al, \be, \ga)=(1,1,2)$, $(2,1,3)$, 
$(2,1,4)$, $(2,1,2)$, $(2,1,4)$ of identity~\eqref{eq:exotique2} are simply
reformulations of Sato's counter-exam\-ples~\eqref{eq:sato1}
up to~\eqref{eq:sato6} and (i) is proved.

For (ii), the idea is to prove that no specialisation of both sides 
of~\eqref{eq:exotique} can follow from the 120 Thomae relations, at least when 
$\al=\beta$. One may note that~\eqref{eq:exotique} cannot formally be a
consequence of any of  Thomae's relations since two (one would be
enough) of its specialisation are not such consequences. But this
does not rule out the possibility that some other
specialisations would follow from Thomae's relations. 
However, we show that this is never 
the case when $\al=\be$ is a positive integer.

We first determine the 120 Thomae relations for the left-hand
side of~\eqref{eq:exotique} when $\al=\be$ and to do this painlessly, we express  
\begin{equation}
\label{eq:membre de gauche}
{} _{3} F _{2} \!\left [ \begin{matrix} {\al+1, \al+1, 2\al}
\\ {3\al+1, 3\al+1}\end{matrix};1
\right ]
\end{equation}
in the symmetric form~\eqref{eq:hardy?} in Theorem~\ref{prop:hardydixonandcie}, which gives
$2x_1=2x_2=2x_4=2x_5=3\al+1$ and $2x_3=4\al.$ The permutations of $x_1$, 
$x_2$, $x_3$, $x_4$ and $x_5$ show that~\eqref{eq:membre de gauche} 
is related only to the 
${} _{3} F_{2}$-series
\begin{equation}
_3F_2\!\!\left[ \begin{matrix} {\al+1, \al+1, 2\al}
\\ {3\al+1, 3\al+1}\end{matrix};1
\right ], \quad
_3F_2\!\!\left [ \begin{matrix} {2\al, 2\al, 2\al}
\\ {3\al+1, 4\al}\end{matrix};1
\right ],\label{eq:fin1}
\end{equation}
and those obtained by the trivial symmetries. 
The same process applied to  the right-hand
side of~\eqref{eq:exotique} (for $\al=\be$),
\begin{equation}
\label{eq:membre de droite}
{} _{3} F _{2} \!\left [ \begin{matrix} {\al, \al, 2\al}
\\ {3\al+1, 3\al+1}\end{matrix};1 
\right ],
\end{equation}
gives $2x_1=2x_2=3\al+2$, $2x_3=4\al+2$ and $2x_4=2x_5=3\al+1$.
The permutations of $x_1$, 
$x_2$, $x_3$, $x_4$ and $x_5$ show that~\eqref{eq:membre de droite}  
is related to
the five  ${} _{3} F _{2}$-series 
\begin{gather}
 _3F_2\!\!\left[ \begin{matrix} {\al, \al, 2\al}
\\ {3\al+1, 3\al+1}\end{matrix};1
\right ],\;
_3F_2\!\!\left [ \begin{matrix} {2\al+2, 2\al+1, 2\al+1}
\\ {4\al+2, 3\al+2}\end{matrix};1
\right ],\;
_3F_2\!\!\left [ \begin{matrix} {2\al+1, \al+1, \al}
\\ {3\al+2, 3\al+1}\end{matrix};1
\right ],\; \nonumber\\
 _3F_2\!\!\left [ \begin{matrix} {2\al+1, 2\al+1, 2\al}
\\ {4\al+2, 3\al+1}\end{matrix};1
\right ],\;
_3F_2\!\!\left [ \begin{matrix} {2\al+2, \al+1, \al+1}
\\ {3\al+2, 3\al+2}\end{matrix};1
\right ], 
\label{eq:fin2}
\end{gather} 
and those obtained by the trivial symmetries. 

Inspection  quickly reveals the impossibility of any numerical coincidence 
between one of the two arrays
of parameters in~\eqref{eq:fin1} and one 
of the five arrays in~\eqref{eq:fin2}, even 
with trivial symmetries. However, each such 
coupling provides a variation 
of~\eqref{eq:exotique} and, for example, we have that
\begin{equation}\label{eq:couplage}
{} _{3} F _{2} \!\left [ \begin{matrix} {2\al, 2\al, 2\al}
\\ {4\al, 3\al+1}\end{matrix};1 
\right ]= \frac{\al(2\al+1)}{(3\al+1)(4\al+1)}\,
{} _{3} F _{2} \!\left [ \begin{matrix} {2\al+2, 2\al+1, 2\al+1}
\\ {4\al+2, 3\al+2}\end{matrix};1 
\right ],
\end{equation}
which will be used below.

We are now in a position to prove the claim about the infinity of 
counter-examples to the cases (i) 
and (ii) of Conjecture~\ref{conj:1}. First, 
thanks to~\eqref{eq:fin3}, we have that 
$$
I(2\al-1,2\al-1,\al, 2\al-1,\al)=
\frac{(2\al-1)!^3\,\al!}{(4\al-1)!\,(3\al)!} \,
{} _{3} F _{2} \!\left [ \begin{matrix} {2\al, 2\al, 2\al}
\\ {4\al, 3\al+1}\end{matrix} ;1
\right ]
$$
and
$$
I(2\al+1,2\al-1,\al, 2\al,\al-1)=
\frac{(2\al+1)!\,(2\al)!\,(2\al-1)!\,\al!}{(4\al+1)!\,(3\al+1)!}\,
{} _{3} F _{2} \!\left [ \begin{matrix} {2\al+2, 2\al+1, 2\al+1}
\\ {4\al+2, 3\al+2}\end{matrix};1
\right ].
$$
We can relate these two equations by~\eqref{eq:couplage}, and the 
simplification of the Gamma-factors yields
$$
I(2\al-1,2\al-1,\al, 2\al-1,\al)= I(2\al+1,2\al-1,\al, 2\al,\al-1),
$$
which is exactly the identity~\eqref{eq:infinicontreexemples} we are looking for.
For both integrals, the non-negativity conditions~\eqref{eq:nonnegative1} and~\eqref{eq:nonnegative2} 
in case (ii) of Conjecture~\ref{conj:1} are verified and 
the above discussion proves that there 
exists no permutation $\rho$ in
the group ${\bf \Phi}$
(and, {\em a fortiori}, also none in ${\bf T}$) such that
$\rho(2\al-1)=2\al+1$, $\rho(2\al-1)=2\al-1$, $\rho(\al)=\al$, $\rho(2\al-1)=2\al$, $\rho(\al)=\al-1.$
Thus, for each value of the positive integer $\al$, we obtain a counter-example to the cases (i) and (ii) of 
Conjecture~\ref{conj:1} at the same time. That this provides
infinitely many counter-examples is a consequence of the 
fact that $I(2\al-1,2\al-1,\al, 2\al-1,\al)$ tends to 0 as $\al$ tends
to infinity, because 
\begin{align*}
\lim_{\al\to +\infty}
I(2\al-1,2\al-1,&\al, 2\al-1,\al)^{1/\al}\\
&=\max_{(x,y)\in [0,1]^2}
\bigg(\frac{x^2(1-x)^2y^2(1-y)}{(1-xy)^2}\bigg)=17-12\sqrt{2}<1.
\end{align*}

\begin{remark} We could do the same thing with $\alpha$ not
necessarily equal to $\beta$.
To find all Thomae relations for the left-hand side
of~\eqref{eq:exotique}, one  
should use Theorem~\ref{prop:hardydixonandcie} with
$$
2x_1=2x_4=2\al+\be+1, \;2x_2=2x_5=2\be+\al+1,\; 2x_3=2\al+2\be,
$$
leading to five different arrays up to trivial symmetries, and for the
right-hand side with
\begin{align*}
2x_1=2\al+\be+2, \;2x_2&=2\be+\al+2,\\ &2x_3
=2\al+2\be+2,\; 2x_4=2\al+\be+1, \;2x_5=2\be+\al+1,
\end{align*}
leading to a complete set of 120 different arrays for 
generic $\al$ and $\beta$ (in fact, only 10 arrays, 
up to trivial symmetries). 
This explains why we consider only the case $\al=\be$, which is 
much simpler to deal with.
\end{remark}

\section{Proof of Theorem~\ref{thm:gamma2}}
\label{sec:contig}

In order to derive Theorem~\ref{thm:gamma2}, we require the following
proposition, relating two ``contiguous" $_3F_2$-series in a way that the
``rest" is a closed form expression. (See the Remark after the proof
of the proposition for an explanation of the term ``contiguous.")

\begin{prop} \label{prop:gamma2}
For any complex numbers $a,b,c$ such that $\Re(d-b-c+1)>0$, and such that 
$a+1$ and $d$ are
not non-positive integers, we have the identity 
\begin{multline} \label{eq:b11} 
  {} _{3} F _{2} \!\left [ \begin{matrix} { a, b, c}\\ { a+1, d}\end{matrix} ;
     {\displaystyle 1}\right ]=
\frac{   ( a - b+1 )  ( a - b+2 )  
    ( a-c+1 )  ( d -1)}
{( a+1 )  (  b-1 )  ( a - d +2 )  
    ( a - d +1 )  }  
    {} _{3} F _{2} \!\left [ \begin{matrix} { a+1,  b-1, c}\\ { a+2, 
     d-1}\end{matrix} ; {\displaystyle 1}\right ]
\\+\frac {(1 - a - a^2 - b + c + a c - b c - d + b d)}
 {  (  b-1 )  ( a - d +2 )  
    ( a - d +1 )}
\frac{ 
     \Gamma({ \textstyle d}) \,\Gamma({ \textstyle d - b - c + 1}) }{
     \Gamma({ \textstyle d-b}) \,\Gamma({ \textstyle  d-c}) }.
\end{multline}
\end{prop}

\begin{proof}
We start by applying the contiguous relation 
\begin{multline} \label{eq:C55}
{} _{3} F _{2} \!\left [ \begin{matrix} { A_1, A_2, A_3}\\ { B_1, B_2}\end{matrix} ;
   {\displaystyle z}\right ]  =
  {\frac{\left( 1 - A_1 + A_2 \right)  \left( B_1 - 1 \right)  
        } 
     {\left( A_1 - 1 \right)  \left( 1 + A_2 - B_1 \right) }}
  {} _{3} F _{2} \!\left [ \begin{matrix} { A_1 - 1, A_2, A_3}\\ { B_1 - 1,
        B_2}\end{matrix} ; {\displaystyle z}\right ] \\+ 
   {\frac {A_2 \left( B_1 - A_1 \right)  
        }
     {\left( A_1 - 1 \right)  \left( B_1 - A_2 - 1  \right) }}
  {} _{3} F _{2} \!\left [ \begin{matrix} { A_1 - 1, A_2 + 1, A_3}\\ { B_1,
        B_2}\end{matrix} ; {\displaystyle z}\right ],
\end{multline}
with $A_1=b$,
$A_2=a$, and $B_1=d$, to the $_3F_2$-series on the left-hand side. 
This yields the expression
$$    \frac{( a-b+1 )  ( d-1 )  
      }{( b-1 )  
      ( a-d+1 ) } \,
   {} _{3} F _{2} \!\left [ \begin{matrix} {  a,b-1, c}\\ { d-1, 
       a+1}\end{matrix} ; {\displaystyle 1}\right ]
-\frac{a ( d-b )  
    }{( b-1 )  
      ( a-d+1 ) }\,
   {} _{2} F _{1} \!\left [ \begin{matrix} { b-1, c}\\ {
       d}\end{matrix} ; {\displaystyle 1}\right ] .
$$
We sum the $_2F_1$-series by means of the Gau{\ss} summation formula
(see \cite[(1.7.6); Appendix (III.3)]{SlatAC}) 
\begin{equation} \label{eq:Gausz} 
{} _{2} F _{1} \!\left [ \begin{matrix} { a, b}\\ { c}\end{matrix} ; {\displaystyle
   1}\right ]  = \frac {\Gamma ( c)\,\Gamma( c-a-b)} {\Gamma( c-a)\, 
  \Gamma( c-b)}.
\end{equation}
Thus, we obtain
$$\frac{( a-b+1 )  ( d-1 )  
     }{( b-1 )  
      ( a-d+1 ) } \, {} _{3} F _{2} \!\left [ \begin{matrix} {  a, b-1,c}\\ { 
        a+1,d-1}\end{matrix} ; 
    {\displaystyle 1}\right ] -
   \frac{a\,\Gamma({ \textstyle d}) \,\Gamma({ \textstyle d-b-c+1}) }
    {( b-1 )  ( a-d+1 ) \,
      \Gamma({ \textstyle d-b}) \,\Gamma({ \textstyle d-c}) }.
$$
Next we apply the contiguous relation
\begin{equation} \label{eq:C15}
{} _3 F _2 \!\left [ \begin{matrix} { A_1, A_2,A_3}\\ { B_1,B_2}\end{matrix} ; {\displaystyle
z}\right ] = {} _3 F _2 \!\left [ \begin{matrix} { A_1 + 1, A_2,A_3}\\ {
    B_1,B_2}\end{matrix} ; {\displaystyle z}\right ]  - 
   {z }
   {\frac{A_2A_3 }
    {B_1B_2}}
   {} _3 F _2 \!\left [ \begin{matrix} { A_1 + 1,A_2+1,A_3+1}\\ {B_1+1,B_2+1}\end{matrix} ;
        {\displaystyle z}\right ],
\end{equation}
with $A_1=a$. This gives
\begin{multline*}
 \frac{( a-b+1 )  ( d-1 )  
      }{( b-1 )  
      ( a-d+1 ) } \, {} _{2} F _{1} \!\left [ \begin{matrix} { b-1, c}\\ { 
       d-1}\end{matrix} ; {\displaystyle 1}\right ]-
   \frac{( a-b+1 )  c 
      }{( a+1 )  
      ( a-d+1 ) }\,
 {} _{3} F _{2} \!\left [ \begin{matrix} { a+1, b, c+1}\\ { a+2,
       d}\end{matrix} ; {\displaystyle 1}\right ] \\+
   \frac{a\,\Gamma({ \textstyle d}) \,\Gamma({ \textstyle d-b-c+1}) }
    {( b-1 )  ( -1 - a + d ) \,
      \Gamma({ \textstyle d-b}) \,\Gamma({ \textstyle d-c}) }.
\end{multline*}
Of course, the $_2F_1$-series can be summed by means of the Gau{\ss}
summation formula \eqref{eq:Gausz}. After some simplification, we
arrive at
\begin{equation} \label{eq:b13a} 
- \frac{( a-b+1 )  c 
      }{( a+1 )
           ( a-d+1 ) }  \,
  {} _{3} F _{2} \!\left [ \begin{matrix} { a+1, b, c+1}\\ { 2 +
         a, d}\end{matrix} ; {\displaystyle 1}\right ]  +
   \frac{( 1 + a + c - d ) \,\Gamma({ \textstyle d}) \,
      \Gamma({ \textstyle d-b-c}) }{( a-d+1 ) \,
      \Gamma({ \textstyle d-b}) \,\Gamma({ \textstyle d-c}) }.
\end{equation}
Now we apply another time the contiguous relation \eqref{eq:C55}, this
time with $A_1=b$, $A_2=a+1$, and $B_1=d$. We obtain
\begin{multline*}
-
   \frac{( a-b+1 )  ( a-b+2 )  c 
      ( d-1 ) }{( a+1 )  ( b-1 )  
      ( a-d+1 )  ( a-d+2 ) } \,
 {} _{3} F _{2} \!\left [ \begin{matrix} { 
       a+1, b-1, c+1}\\ { a+2, d-1}\end{matrix} ; {\displaystyle
       1}\right ] \\+  \frac{( a-b+1 )  c ( d-b )  
        }{( b-1 )  
        ( a-d+1 )  ( a-d+2 ) } \,
 {} _{2} F _{1} \!\left [ \begin{matrix} { b-1, c+1}\\ {
         d}\end{matrix} ; {\displaystyle 1}\right ]  \\+
   \frac{( a + c - d +1) \,\Gamma({ \textstyle d}) \,
      \Gamma({ \textstyle d-b-c}) }{( a-d+1 ) \,
      \Gamma({ \textstyle d-b}) \,\Gamma({ \textstyle d-c}) },
\end{multline*}
and after evaluation of the $_2F_1$-series by means of Gau{\ss}'
summation formula \eqref{eq:Gausz},
\begin{multline*}
 - \frac{( a-b+1 )  ( a-b+2 )  c 
        ( d-1 )  
      }{
        ( a+1 )  ( b-1 )  
        ( a-d+1 )  ( a-d+2 ) }  \,
  {} _{3} F _{2} \!\left [ \begin{matrix} { a+1, b-1, c+1}\\ {
         a+2, d-1}\end{matrix} ; {\displaystyle 1}\right ] \\ +
   \frac{P(a,b,c,d)\,
     \Gamma({ \textstyle d}) \,\Gamma({ \textstyle d-b-c}) }
      {( b-1 )  ( a-d+1 )  
      ( a-d+2 ) \,\Gamma({ \textstyle d-b}) \,
      \Gamma({ \textstyle d-c}) },
\end{multline*}
where 
\begin{multline*}
P(a,b,c,d)=-2 - 3 a - a^2 + 2 b + 3 a b + a^2 b - 3 c - 2 a c +
        3 b c + a b c - c^2 - a c^2 \\+ b c^2 + 3 d + 2 a d -
        3 b d - 2 a b d + 2 c d + a c d - 2 b c d - d^2 + b d^2.
\end{multline*}
The final contiguous relation that we apply is
\begin{multline} \label{eq:C27}
{} _3 F _2 \!\left [ \begin{matrix} { A_1, A_2, A_3}\\ { B_1,B_2}\end{matrix} ;
   {\displaystyle z}\right ] = 
  {\frac{ A_1 - A_2 - 1  
        } {A_1 - 1}}\,
  {} _3 F _2 \!\left [ \begin{matrix} { A_1 - 1, A_2, A_3}\\ { B_1,B_2}\end{matrix} ;
        {\displaystyle z}\right ] \\+ 
   {\frac{A_2  } {A_1 - 1}}\,
  {} _3 F _2 \!\left [ \begin{matrix} { A_1 - 1, A_2 + 1, A_3}\\ {
        B_1,B_2}\end{matrix} ; {\displaystyle z}\right ],
\end{multline}
with $A_1=c+1$ and $A_2=a+1$. The result is
\begin{multline*}
   \frac{( a-b+1 )  ( a-b+2 )  
      ( a-c+1 )  ( d-1 )  
   }{( a+1 )
         ( b-1 )  ( a-d+1 )  
      ( a-d+2 ) } \,
  {} _{3} F _{2} \!\left [ \begin{matrix} { a+1, b-1, c}\\ { a+2,
       d-1}\end{matrix} ; {\displaystyle 1}\right ]  \\
- \frac{( a-b+1 )  ( a-b+2 )  
        ( d-1 )  
        }{( b-1 )  
        ( a-d+1 )  ( a-d+2 ) }  \,
{} _{2} F _{1} \!\left [ \begin{matrix} { b-1, c}\\ { 
         d-1}\end{matrix} ; {\displaystyle 1}\right ]  
\\+
   \frac{P(a,b,c,d)\,\Gamma({ \textstyle d}) \,\Gamma({ \textstyle d-b-c}) }
      {( b-1 )  ( a-d+1 )  
      ( a-d+2 ) \,\Gamma({ \textstyle d-b}) \,
      \Gamma({ \textstyle d-c}) }.
\end{multline*}
A last use of Gau{\ss}' summation formula and some simplification then
leads to \eqref{eq:b11}.
\end{proof}

\begin{remark}
Since we shall re-use it in
Section~\ref{sec:contig2}, it will be beneficial if we briefly summarise
the idea of the proof of the above proposition: it is crucially based
on the fact that the $_3F_2$-series on the left of \eqref{eq:b11} has
the parameter $a$ on top and the parameter $a+1$ at the bottom. Now we
apply elementary contiguous relations (such as the one in
\eqref{eq:C55}). In principle, it expresses our $_3F_2$-series as a
sum of two other $_3F_2$-series in which the parameters are
``contiguous" to the original $_3F_2$-series, meaning that they differ
from the parameters of the original series by small integer amounts.
(In \eqref{eq:C55}, these differences are $0$ and $\pm1$.) However, in
one of the two $_3F_2$-series on the right-hand side of the relation,
the top parameter $A_2=a$ is raised by $1$, while the bottom parameter
$B_2=a+1$ is left invariant. Thus, the two $(a+1)$'s cancel, and the
$_3F_2$-series reduces to a $_2F_1$-series, to which the Gau{\ss}
summation formula \eqref{eq:Gausz} can be applied to express it in
closed form. This partial simplification happens as well when we apply
the contiguous relations \eqref{eq:C15} and \eqref{eq:C27}. Thus, each
time, we obtain a $_3F_2$-series plus an additional expression in closed
form. These additional expressions are put together, and they finally
form the expression containing the gamma functions on the right-hand
side of \eqref{eq:b11}. However, since several similar, but not
identical, such expressions were put together, 
when factoring the resulting term, 
a polynomial factor built up. Hence, in order to obtain a relation
between two $_3F_2$-series without any additional term,  
this polynomial factor must vanish. While, normally 
(i.e., if one plays
the above described game in a random fashion), equating this
polynomial factor to zero will not have any nice solutions 
(in
particular, no {\it integral\/} solutions, which we would however need
to construct counter-examples to the conjecture by Rhin and Viola, in
the cases of Propositions~\ref{prop:gamma2}--\ref{prop:gamma3}), the
contiguous relations have been carefully selected so that at least
one of the variables $a,b,c,d$ is contained only linearly in the
polynomial factor. This makes it possible to have many non-trivial
solutions when equating the polynomial factor to zero.
\end{remark}

In view of the above remark,
the proof of Theorem~\ref{thm:gamma2} is now straight-forward.

\begin{proof}[Proof of Theorem~\ref{thm:gamma2}]
If we now choose $d$ such that the polynomial factor on the right-hand
side of \eqref{eq:b11} vanishes, that is,
$$
d= c+\frac{a(  a - c+1 ) }{b-1} +1,
$$
and subsequently do the replacements $a\to\al$, $b\to\be$, $c\to\ga$,
then we obtain exactly \eqref{eq:b12}.
\end{proof}

\section{Proof of Theorem~\ref{theo:2a}} \label{sec:proofa}

As already mentioned in Section~\ref{sec:pattern}, the case 
$(\al, \be, \ga)=(4,3,4)$ of identity~\eqref{eq:b12} is a simple
reformulation of Sato's counter-exam\-ple~\eqref{eq:sato3}.
Thus, (i) is proved.

For (ii), we proceed in a similar fashion as in
Section~\ref{sec:proof}. First of all, we observe that identity
\eqref{eq:infinicontreexemplesa} is the special case of \eqref{eq:b12}
in which $\al$ is replaced by $\al^2$, and in which
$\be=\al+1$ and $\ga=\al^2$, again via the translation
\eqref{eq:lossofsymmetry}. To wit, this is
\begin{equation} \label{eq:b3} 
 {} _{3} F _{2} \!\left [ \begin{matrix} { \al^2,\al^2,
   \al+1}\\ { \al ^2+1, 
   \al^2+\al +1}\end{matrix} ; {\displaystyle 1}\right ]
=
 \frac{\al ^3+1}{\al ^2+1} \;
  {} _{3} F _{2} \!\left [ \begin{matrix} {  \al^2 
      +1, \al^2, \al}\\ { 
      \al^2+2, \al^2+\al}\end{matrix} ; {\displaystyle 1}\right ] .
\end{equation}
In the sequel we concentrate on this special case, always assuming
that $\al$ is a positive integer strictly greater than $1$.

Using Thomae's relations, we can generate three 
other series which are related to the $_3F_2$-series on the left-hand
side of \eqref{eq:b3}, namely
\begin{equation} \label{eq:rel1}
 {} _{3} F _{2} \!\left [ \begin{matrix}
\al+1,\al+1,\al+1\\\al+2, \al^2+\al+1
\end{matrix}; {\displaystyle 1}\right ] ,\quad 
 {} _{3} F _{2} \!\left [ \begin{matrix}
1,1,\al+1\\\al+2,\al^2+1
\end{matrix}; {\displaystyle 1}\right ] ,\quad 
 {} _{3} F _{2} \!\left [ \begin{matrix}
1, \al^2, \al^2-\al\\\al^2+1, \al^2+1
\end{matrix}; {\displaystyle 1}\right ] .
\end{equation}
On the other hand, there are six series related to the $_3F_2$-series on
the right-hand side of \eqref{eq:b3},
\begin{multline} \label{eq:rel2}
 {} _{3} F _{2} \!\left [ \begin{matrix}
\al-1,\al^2,\al^2\\\al^2+1,\al^2+\al
\end{matrix}; {\displaystyle 1}\right ] ,\quad 
 {} _{3} F _{2} \!\left [ \begin{matrix}
\al-1,\al,\al\\\al+1,\al^2+\al
\end{matrix}; {\displaystyle 1}\right ] ,\quad 
 {} _{3} F _{2} \!\left [ \begin{matrix}
2,\al^2+1,\al^2-\al+2\\\al^2+2,\al^2+2
\end{matrix}; {\displaystyle 1}\right ] ,\\
 {} _{3} F _{2} \!\left [ \begin{matrix}
1,\al^2,\al^2-\al+2\\\al^2+1,\al^2+2
\end{matrix}; {\displaystyle 1}\right ] ,\quad  
 {} _{3} F _{2} \!\left [ \begin{matrix}
1,2,\al\\\al+1,\al^2+2
\end{matrix}; {\displaystyle 1}\right ] ,\quad  
 {} _{3} F _{2} \!\left [ \begin{matrix}
1,1,\al-1\\\al+1,\al^2+1
\end{matrix}; {\displaystyle 1}\right ] .
\end{multline}
None of these match with the series on the left-hand side of
\eqref{eq:b3} or with one of the series in \eqref{eq:rel1}. Thus,
indeed, for any positive integer $\al$,
\eqref{eq:infinicontreexemplesa} is a counter-example to the
conjecture by Rhin and Viola. 

Finally, in order to see that \eqref{eq:infinicontreexemplesa} 
produces infinitely many
counter-examples, we show again that the involved integral tends to
zero when $\al$ tends to infinity. Indeed, for $\al\ge 1$, we have 
\begin{align*} \notag
I(\al^2-1,\al-1,\al^2-\al+1,\al-1,0)
&=\int_0^1 \!\!\int_0^1
\frac{x^{\al^2-1}(1-x)^{\al-1}y^{\al-1}(1-y)^{\al^2-\al+1}}
{(1-xy)^{\al^2+1}}\,\dd x\, \dd y\\
\notag
&\le \int_0^1 \!\!\int_0^1 x^{\al -1}y ^{\al-1} \frac{\dd x\, \dd y}{1-xy} = 
\sum_{k=0} ^{\infty} \frac{1}{(k+\al)^2},
\end{align*}
from which the claim follows. (In the second line, we used 
the trivial facts that $x^{\al^2}\le x ^\al$ and 
$(1-x)^{\al-1}(1-y)^{\al ^2-\al+1}\le (1-xy)^{\al^2}$ for $0\le x, y\le 1$.)
This completes the proof of Theorem~\ref{theo:2a}.

\begin{remark}
It is obvious that Theorem~\ref{thm:gamma2} will generate many more
counter-examples to the conjecture by Rhin and Viola, by choosing the
parameters $\al,\be,\ga$ to be positive integers in other ways
such that $\al ( \al - \ga +1)/({\be-1}+1 )$ is as well a positive
integer (and such that the conditions
\eqref{eq:cond1}--\eqref{eq:cond2} are satisfied). To have a
convenient parametrisation, one would replace $\ga$ by $\al+1-\ga$,
subsequently $\al$ by $a_1a_2$, $\ga$ by
$c_1c_2$, and $\be$ by $a_1c_1+1$. The resulting relation is
\begin{multline} 
\label{eq:b14}
 {} _{3} F _{2} \!\left [ \begin{matrix} { a_1a_2, 
    a_1c_1+1,  a_1a_2 -
    c_1 c_2+1}\\ {  a_1a_2+1, 
    a_1a_2 + a_2 c_2 -
   c_1 c_2+2}\end{matrix} ; {\displaystyle 1}\right ] \\=
   \frac{(  a_1a_2 -
        a_1c_1 +1)  
     (  a_1a_2 +
        a_2 c_2 - c_1 c_2
       +1 ) }
     {( a_1a_2+1 )  
      (  a_2 c_2 -
        c_1 c_2+1 ) }
\,{} _{3} F _{2} \!\left [ \begin{matrix} { 
       a_1a_2+1, a_1c_1, 
      a_1a_2 - c_1 c_2+1}\\ { 
       a_1a_2+2,  a_1a_2 +
       a_2 c_2 -
      c_1 c_2+1}\end{matrix} ; {\displaystyle 1}\right ] .
\end{multline}
\end{remark}

\section{More exotic contiguous relations}
\label{sec:contig2}

In this section, we present two more relations of the kind of
Theorem~\ref{thm:gamma2} (which itself followed from the more
general Proposition~\ref{prop:gamma2}), see Theorems~\ref{thm:beta} and
\ref{thm:gamma1}. These are obtained along the
lines described in the Remark after the proof of
Proposition~\ref{prop:gamma2}. The two theorems 
imply further counter-examples to the conjecture by Rhin and Viola.

\begin{prop} \label{prop:beta}
For any complex numbers $a,b,c$ such that $\Re(d-b-c+1)>0$, and such that 
$a+1$ and $d$ are
not non-positive integers, we have the identity 
\begin{multline} \label{eq:b5} 
 {} _{3} F _{2} \!\left [ \begin{matrix} { a, b, c}\\ { a+1, 
    d}\end{matrix} ; {\displaystyle 1}\right ]=
  \frac{b c \left( a-d-1 \right)  \left( a - d \right)  
     }{\left( a-b
           \right)  \left( a - c \right)  d \left( d+1 \right) } 
      {} _{3} F _{2} \!\left [ \begin{matrix} { a, b+1, c+1}\\ { 
         a+1, d+2}\end{matrix} ; {\displaystyle 1}\right ]\\
     +\frac {a ( b c + a d - b d - c d ) } 
  {      ( a - b )  ( a - c ) }
 \frac{
      \Gamma({ \textstyle d}) \,\Gamma({ \textstyle d-b-c+1}) }{
      \Gamma({ \textstyle d-b+1}) \,\Gamma({ \textstyle d-c+1}) }.
\end{multline}
\end{prop}
\begin{proof}To the left-hand side, we apply the contiguous relation
\begin{multline} \label{eq:C54}
{} _3 F _2 \!\left [ \begin{matrix} { A_1, A_2, A_3}\\ { B_1, B_2}\end{matrix} ;
   {\displaystyle z}\right ] =
  {\frac {A_2 \left( B_1 - A_1 \right)   }
     {\left( A_2 - A_1 \right)  B_1}}
  {} _3 F _2 \!\left [ \begin{matrix} { A_1, A_2 + 1,
        A_3}\\ { B_1 + 1, B_2}\end{matrix} ; {\displaystyle z}\right ] \\+ 
   {\frac {A_1 \left( B_1 - A_2 \right)  
        }
     {\left( A_1 - A_2 \right)  B_1}}
  {} _3 F _2 \!\left [ \begin{matrix} { A_1 + 1, A_2, A_3}\\ { B_1 + 1,
        B_2}\end{matrix} ; {\displaystyle z}\right ]
\end{multline}
with $A_1=b$, $A_2=c$, and $B_1=d$. As a result we obtain
\begin{equation*}
 \frac{c( d-b ) 
      }{( c-b ) d}
  {} _{3} F _{2} \!\left [ \begin{matrix} { a, b, c+1}\\ { a+1, 
       d+1}\end{matrix} ; {\displaystyle 1}\right ]
   + \frac{b ( d-c )  
    }{( b - c)d } 
 {} _{3} F _{2} \!\left [ \begin{matrix} { a, b+1, c}\\ { a+1, 
       d+1}\end{matrix} ; {\displaystyle 1}\right ] .
\end{equation*}
We apply the contiguous relation \eqref{eq:C54} again, to the first
series with $A_1=a$, $A_2=b$, and $B_1=d+1$, to the second
 with $A_1=a$, $A_2=c$, and $B_1=d+1$. After some simplification, this
leads to the expression
\begin{multline*}
-
   \frac{b c \left( a-d-1 \right)  \left( a - d \right)  
     }{\left( -a + b \right)
         \left( a - c \right)  d \left( d+1 \right) }
 {} _{3} F _{2} \!\left [ \begin{matrix} { a, b+1, c+1}\\ { a+1,
       d+2}\end{matrix} ; {\displaystyle 1}\right ]\\
+ \frac{a c \left( d-b \right)  \left( d-b+1 \right)  
      }{\left( a - b \right)  
      \left( c-b \right)  d \left( d+1 \right) } 
 {} _{2} F _{1} \!\left [ \begin{matrix} { b, c+1}\\ { 
       d+2}\end{matrix} ; {\displaystyle 1}\right ]+
   \frac{a b \left( d-c \right)  \left( d-c+1 \right)  
      }{\left( a - c \right)  
      \left( b - c \right)  d \left( d+1 \right) } 
 {} _{2} F _{1} \!\left [ \begin{matrix} { b+1, c}\\ { 
       d+2}\end{matrix} ; {\displaystyle 1}\right ] .
\end{multline*}
Finally, we use Gau{\ss}' summation formula \eqref{eq:Gausz} to
evaluate the two $_2F_1$-series. Some manipulation then yields the
claimed result on the right-hand side of \eqref{eq:b5}.
\end{proof}
We may now choose $a$ so that the second term on the right-hand side
of \eqref{eq:b5} vanishes, that is, we choose
$$a=b+c-\frac {b c}d.$$
After the additional replacements of $b$ by $\be$, of $c$ by $\ga$,
and of $d$ by $\be\ga/\de$, we arrive at the
following result.

\begin{theo} \label{thm:beta}
For any complex numbers $\al,\be,\ga$ such that  
$\Re\(\frac {\be\ga} {\de}-\be-\ga+1\)>0$, and such that
$\be + \ga - \de+1$ and $\frac {\be\ga} {\de}$ are
not non-positive integers, we have the identity 
\begin{equation} \label{eq:b6} 
 {} _{3} F _{2} \!\left [ \begin{matrix} { \be + \ga - \de, \be,
\ga}\\ {  \be + \ga -
    \de+1, \frac{\be \ga}{\de}}\end{matrix} ; {\displaystyle 1}\right ]  =
   \frac{ \be \ga + \de - \be \de - \ga \de + \de^2 
   }{\be \ga + \de} \,
  {} _{3} F _{2} \!\left [ \begin{matrix} { \be + \ga - \de, \be+1, \ga+1}\\ {
        \be + \ga - \de+1, \frac{\be \ga}{\de}+2}\end{matrix} ; {\displaystyle
       1}\right ] .
\end{equation}
\end{theo}

Again, if one wants a more convenient parametrisation for generating
counter-examples to the conjecture by Rhin and Viola, then one would
replace $\be$ by $b_1b_2$, $\ga$ by $c_1c_2$, and $\de$ by
$b_1c_1$. The resulting relation then is
\begin{multline} 
\label{eq:b15} 
 {} _{3} F _{2} \!\left [ \begin{matrix} { b_1b_2 + c_1c_2 
- b_1c_1, b_1b_2, c_1c_2}\\ {  b_1b_2 + c_1c_2 -
    b_1c_1+1, {b_2 c_2}}\end{matrix}
; {\displaystyle 1}\right ]  \\=
   \frac{ b_1b_2 c_1c_2 + b_1c_1 - 
b_1b_2 b_1c_1 - c_1c_2 b_1c_1 + b_1^2c_1^2 
   }{b_1b_2 c_1c_2 + b_1c_1} \\
\times
  {} _{3} F _{2} \!\left [ \begin{matrix} { b_1b_2 + 
c_1c_2 - b_1c_1, b_1b_2+1, c_1c_2+1}\\ {
       b_1b_2 + c_1c_2 - b_1c_1+1,
{b_2 c_2}+2}\end{matrix} ; {\displaystyle
       1}\right ] .
\end{multline}

%
%

\begin{prop} \label{prop:gamma3}
For any complex numbers $a,b,c$ such that $\Re(d-b-c+1)>0$, and such that 
$a+1$ and $d$ are
not non-positive integers, we have the identity 
\begin{multline} \label{eq:b13} 
{} _{3} F _{2} \!\left [ \begin{matrix} { a, b, c}\\ { a+1,
   d}\end{matrix} ; {\displaystyle 1}\right ]
=
  \frac{( a - b +1)  (a-c+1 )  
        }{( a+1 )  
        ( a - d +1) }  \,
 {} _{3} F _{2} \!\left [ \begin{matrix} { a+1, b,c}\\ { a+2,
         d}\end{matrix} ; {\displaystyle 1}\right ]\\ -
   \frac{\Gamma({ \textstyle d}) \,\Gamma({ \textstyle d - b - c + 1}) }
    {( a - d +1) \,\Gamma({ \textstyle d-b}) \,
      \Gamma({ \textstyle d-c}) }.
\end{multline}
\end{prop}
\begin{proof}
The first few steps of this proof are identical with the one of 
Proposition~\ref{prop:gamma2}. More precisely, we use that the series
on the left-hand side is equal to the expression \eqref{eq:b13a}.
There, we apply now instead the contiguous relation \eqref{eq:C27}
with $A_1=c+1$ and $A_2=a+1$. As a result, we obtain
\begin{multline*}
   \frac{\left( a - b+1 \right) \,\left(a-c+1 \right) \,
     }{\left( a+1 \right) \,
      \left( a - d +1\right) } \,
  {} _{3} F _{2} \!\left [ \begin{matrix} {  a+1, b,c}\\ {  a+2,
       d}\end{matrix} ; {\displaystyle 1}\right ]
-  \frac{ a - b +1
     }{a - d+1} \,
    {} _{2} F _{1} \!\left [ \begin{matrix} {b, c}\\ { d}\end{matrix} ;
         {\displaystyle 1}\right ]\\+
   \frac{\left(  a + c - d+1 \right) \,\Gamma({ \textstyle d}) \,
      \Gamma({ \textstyle  d-b-c}) }{\left( a - d+1 \right) \,
      \Gamma({ \textstyle d-b}) \,\Gamma({ \textstyle  d-c}) },
\end{multline*}
which, by another use of the Gau{\ss} summation formula
\eqref{eq:Gausz} and some
simplification, turns out to be equal to the right-hand side of
\eqref{eq:b13}. 
\end{proof}

An iterative use of Proposition~\ref{prop:gamma3} produces the following
formula.

\begin{coro} \label{prop:gamma1}
For any complex numbers $a,b,c$ such that $\Re(d-b-c+1)>0$, and such that 
$a+1$ and $d$ are
not non-positive integers, we have the identity 
\begin{multline} \label{eq:b9} 
      {} _{3} F _{2} \!\left [ \begin{matrix} { a, b, c}\\ { a+1,
       d}\end{matrix} ; {\displaystyle 1}\right ]      
=
\frac{ ( a - b +1) ( a - b+2 )
    ( a - c+1 ) ( a - c+2 )}
 {( a+1 ) ( a+2 ) 
      ( a-d+2)(a-d+1 ) }\,
    {} _{3} F _{2} \!\left [ \begin{matrix} { a+2, b, c}\\ { a+3,
     d}\end{matrix} ; {\displaystyle 1}\right ]\\
-\frac {(3 + 5 a + 2 a^2 - b - a b - c - a c + b c - d - a d)} 
 {( a+1 ) 
      ( a-d+2)(a-d+1 )}
 \frac{ \Gamma({ \textstyle d}) \,
     \Gamma({ \textstyle d - b - c + 1}) }{
 \Gamma({ \textstyle d-b}) \,
     \Gamma({ \textstyle d-c}) }.
\end{multline}
\end{coro}

If we now choose $d$ such that the polynomial factor on the right-hand
side of \eqref{eq:b9} vanishes, that is,
$$d=2 a - b - c +3+ \frac{b c}{ a+1},$$
then we obtain the following theorem.

\begin{theo} \label{thm:gamma1}
For any complex numbers $\al,\be,\ga$ such that  
$\Re\(2\al-2\be-2\ga+\frac {\be} {\al+1}+4\)>0$, and such that
$\al+1$ and $  2 \al - \be - \ga+ \frac{\be \ga}{\al+1}+3$ are
not non-positive integers, we have the identity 
\begin{multline} \label{eq:b10} 
{} _{3} F _{2} \!\left [ \begin{matrix} { \al, \be, \ga}\\ { \al+1, 
         2 \al - \be - \ga + \frac{\be \ga}{\al+1}+3}\end{matrix} ; {\displaystyle
         1}\right ] \\=
\frac {(\al+1)( \al - \be +2)  ( \al - \ga+2 ) } {( \al +2)  
        (  3 \al + \al^2 - \be - \al \be - \ga - \al \ga + \be \ga+2 )  } 
    {} _{3} F _{2} \!\left [ \begin{matrix} { \al+2, \be, \ga}\\ { \al+3, 
     2 \al - \be - \ga+ \frac{\be \ga}{\al+1}+3}\end{matrix} ; {\displaystyle 1}\right
     ].
\end{multline}
\end{theo}

If one wants a more convenient parametrisation for generating
counter-examples to the conjecture by Rhin and Viola, then one would
replace $\be$ by $b_1b_2$, $\ga$ by $c_1c_2$, and $\al$ by
$b_1c_1-1$. The resulting relation then is
\begin{multline} 
\label{eq:b16}
 {} _{3} F _{2} \!\left [ \begin{matrix} { 
   b_1c_1-1,b_1b_2,
   c_1c_2}\\ {b_1c_1,  2b_1c_1  -
   b_1b_2  
    +b_2c_2 -
   c_1c_2+1}\end{matrix} ; {\displaystyle 1}\right ] \\ =
    \frac{(  b_1c_1-b_1b_2 
       +1 )  
      ( b_1c_1 -
       c_1c_2+1 )  
       }{( b_1c_1+1 )  
      ( b_1c_1 -b_1b_2 
        +b_2c_2 -
       c_1c_2+1 ) } \\
\times {} _{3} F _{2} \!\left [ \begin{matrix} { 
      b_1c_1+1,b_1b_2,
      c_1c_2}\\ { b_1c_1+2,  2b_1c_1
        -b_1b_2   
       +
      b_2c_2 -
      c_1c_2+1}\end{matrix} ; {\displaystyle 1}\right ].
\end{multline}

\section{Postlude: how were these identities found?}
\label{sec:post}

The reader may wonder how we found the identities in
Theorem~\ref{theo:??} and
Propositions~\ref{prop:gamma2}--\ref{prop:gamma3} 
(the latter implying 
Theorems~\ref{thm:gamma2}, \ref{thm:beta} and \ref{thm:gamma1})
and their proofs. This section describes some of the ideas that
led us to their discovery, with some of them being interesting in their own
right, as we believe. Since we shall make reference to it several
times, we mention right away that all the hypergeometric
calculations were carried out using the first author's {\sl
Mathematica} package HYP \cite{hyp}.

The counter-examples \eqref{eq:sato1}--\eqref{eq:sato3} of Sato, in
their original form, do not give any hints for a general result that
may be behind them. However, as we explain in
Section~\ref{sec:pattern}, if we bring them into different, but
equivalent, forms using Thomae's relations, patterns emerge. More
precisely, by staring at the forms \eqref{eq:counter1} and
\eqref{eq:counter2} of \eqref{eq:sato1}--\eqref{eq:sato6}, we
extracted the wild guess that \eqref{eq:b12} should hold.
The first proof that we found (which is not presented here) showed
first the special case $\ga=\al+\be$ of \eqref{eq:b12}, given in 
\eqref{eq:exotique}, by using elementary contiguous relations.
A somewhat involved analytic continuation argument,  
using the Gosper--Zeilberger algorithm (see below) and Carlson's theorem 
then extended \eqref{eq:exotique} to \eqref{eq:exotique2}. 

However, it was
``obvious" to us that one should be able to prove \eqref{eq:b12} by a
combination of several classical transformation formulae for
hypergeometric series. Clearly, since we know that \eqref{eq:b12} is
not a consequence of Thomae's relations, the classical
$_3F_2$-transformations are not of any use. So we asked HYP to tell us
which (of the built-in) transformations can be applied to the
left-hand side of \eqref{eq:b12}. (This is done by using {\tt TListe};
see \cite{hyp}.) The only ``non-standard" transformation that HYP came
up with was \eqref{eq:T3240}. (This is {\tt T3240} in HYP.) So 
we applied it and quickly realized that we could exchange $\al$ and
$\al+1$ in the obtained $_7F_6$-series (cf.\ \eqref{eq:7F6a} and
\eqref{eq:7F6b}) and apply \eqref{eq:T3240} in the other direction, in
order to obtain a result different from the original $_3F_2$-series,
which then turned out to be exactly the right-hand side of \eqref{eq:b12}.

Having found an explanation for $83.33333\dots$ percent of Sato's
counter-examples did not completely satisfy us. We also wanted an
explanation for \eqref{eq:sato3}. Since this is just one single
identity, there is only very little guidance where to look for. 
What caught
our eyes was that, in the hypergeometric form \eqref{eq:sato3hyp},
both $_3F_2$-series were balanced (that is, the sum of the lower
parameters exceeds the sum of the upper parameters by exactly 1). Not
only that, in both series there is a lower parameter which exceeds an
upper parameter by exactly 1. So, we made our computer work out the
values of all series of the form 
$${} _{3} F _{2} \!\left [ \begin{matrix} { a, b, c}\\ { a+1,
      b+c}\end{matrix} ; {\displaystyle 1}\right ]$$
for $1\le a,b,c\le 40$,
and then compared which series were rational multiples of each other.
By staring at the results, we extracted identities such as
\begin{equation} 
\label{eq:id1} 
 {} _{3} F _{2} \!\left [ \begin{matrix} { \al^2,
   \al+1,\al^2}\\ { \al ^2+1, 
   \al^2+\al +1}\end{matrix} ; {\displaystyle 1}\right ]
=
\frac
{\al^3+1}
{\al^2 +1}\,
  {} _{3} F _{2} \!\left [ \begin{matrix} {  \al^2 
      +1, \al, \al^2}\\ { 
      \al^2+2, \al^2+\al}\end{matrix} ; {\displaystyle 1}\right ] 
\end{equation}
(this is identity \eqref{eq:b3}, the special case $\al\to\al^2$,
$\be=\al+1$, $\ga=\al^2$ of \eqref{eq:b12}), or
\begin{equation} 
\label{eq:id2} 
 {} _{3} F _{2} \!\left [ \begin{matrix} { 
    {\al}^ 2-\al+1, \al, 
    {\al}^ 2-\al}\\ { 
    {\al}^ 2-\al+2, {\al}^ 2}\end{matrix} ;
    {\displaystyle 1}\right ]  =
   \frac{\al  }{{\al}^2+1}\,
{} _{3} F _{2} \!\left [ \begin{matrix} { 
       {\al}^ 2-\al+1, \al+1,  {\al}^ 2-\al+1}\\ { 
       {\al}^ 2-\al+2,  {\al}^2+2
       }\end{matrix} ; {\displaystyle 1}\right ]
\end{equation}
(this is the special case $\be=\al$,
$\ga=\al^ 2-\al$, $\de=\al-1$ of \eqref{eq:b6}) or
\begin{equation} \label{eq:id3} 
 {} _{3} F _{2} \!\left [ \begin{matrix} 
    {6\al+1, 4\al+2,3\al+1}\\ { 
    6\al+2,7\al+3}\end{matrix} ;
    {\displaystyle 1}\right ]  =
   \frac{3\al+2  }{3{\al}+3}\,
{} _{3} F _{2} \!\left [ \begin{matrix} 
    {6\al+3, 4\al+2,3\al+1}\\ { 
    6\al+4,7\al+3}\end{matrix} ;
    {\displaystyle 1}\right ]
\end{equation}
(this is the special case $\al\to6\al+1$, $\be=4\al+2$, $\ga=3\al+1$
of \eqref{eq:b10}).

We then attempted to prove these identities. It seems sort of
``obvious" that one should be able to prove them by using known contiguous
relations. Indeed, in HYP there are approximately 100 such
contiguous relations built-in. We played with those, but  
we were not able to arrive at the right-hand sides of the conjectured
identities. At some point, we had the idea to ``cheat" and to
make recourse to the ``modern" way of treating hypergeometric series,
namely applying the Gosper--Zeilberger algorithm 
(see \cite{ek,PeWZAA,ZeilAP,ZeilAM,ZeilAV}; what we do below is in the
spirit of \cite{PaulAZ}). For example, aiming to
prove (a generalisation of) \eqref{eq:id1}, we considered
the series
\begin{equation} \label{eq:abc1} 
{} _{3} F _{2} \!\left [ \begin{matrix} 
    a+n,b-n,c\\ a+n+1,b+c-n
    \end{matrix} ;
    {\displaystyle 1}\right ]
\end{equation}
and tried to find a first-order recurrence for it (which is what
\eqref{eq:b12} is). Thus, we put the summand of this
series,
$$F(n,k)=\frac {(a+n)_k\,(b-n)_k\,(c)_k} {(a+n+1)_k\,(b+c-n)_k\,k!}$$
into the Gosper--Zeilberger algorithm, and we got
\begin{multline} \label{eq:rek1}
( a + n+1) (b-n-1) (a - b - c + 2 n+1)
         (a - b - c + 2 n+2) F(n,k) \\+
      (a-b+2n+2) (a-b+2n+1) (a-c+n+1)
       ( b + c - n-1) F(n+1,k) =\De_k F(n,k) R(n,k),
\end{multline}
where 
\begin{multline*}
R(n,k)=  \frac {k (b + c + k - n-1 )} {( a + n
       )  (  b + c - n-1 )  (  b + k - n-1 )}
\\\times
      ( k ( b + c - n-1 )  
          (  a + n +1)  
          ( 1 - a - a^2 - 2 b + b^2 + a c + n - 2 a n - 2 b n +
            c n   )\\  +\text{terms not containing $k$}),
\end{multline*}
and where $\De_k$ is the forward difference operator, $(\De_k
f)(k)=f(k+1)-f(k)$. If we now sum both sides of \eqref{eq:rek1} over
$k$ from $0$ to $N$, then we obtain
\begin{multline*}
( a + n+1) (b-n-1) (a - b - c + 2 n+1)
         (a - b - c + 2 n+2)
\sum _{k=0} ^{N} F(n,k) \\+
      (a-b+2n+2) (a-b+2n+1) (a-c+n+1)
       ( b + c - n-1)
\sum _{k=0} ^{N} F(n+1,k)\\ = F(n,N+1) R(n,N+1),
\end{multline*}
since the terms on the right-hand side telescope. Subsequently, the
limit $N\to\infty$ yields
\begin{multline} \label{eq:m=0}  
( a+n+1 )  (  b-n-1 )  ( a - b-c+2n +2 )  
    ( a - b-c +2n+1 )\\
\times
  {} _{3} F _{2} \!\left [ \begin{matrix} { a+n, b-n, c}\\ { a+n+1, 
         b+c-n}\end{matrix} ;
     {\displaystyle 1}\right ]\\=
{   ( a - b+2n+1 )  ( a - b+2n+2 )  
    ( a-c+n+1 )  ( b+c-n -1)}
  \\
\times
    {} _{3} F _{2} \!\left [ \begin{matrix} { a+n+1,  b-n-1, c}\\ { a+n+2, 
     b+c-n-1}\end{matrix} ; {\displaystyle 1}\right ]
\\+( a+n+1 )
(1 - a - a^2 - b + c + a c - b c - d + b d + n - 2 a n - b n +
   2 c n - d n )
\frac{ 
     \Gamma({ \textstyle b+c-n})  }{
     \Gamma({ \textstyle c}) \,\Gamma({ \textstyle  b-n}) }.
\end{multline}
(The reader should notice that this is \eqref{eq:b11} with $a$
replaced by $a+n$, $b$ replaced by $b-n$, and $d$ replaced by
$b+c-n$.)

At this point, we became greedy. Why should this be something special
for balanced series? So, we replaced the bottom parameter $b+c-n$ in
\eqref{eq:abc1} by $d-n$, --- and we were disappointed to learn that
the Gosper--Zeilberger algorithm is unable to find a two-term
recurrence for this more general series. (It finds only a three-term
recurrence.) However, it {\it does} find a
two-term recurrence for every $d$ of the form $d=b+c+m$, where $m$
is a non-negative integer. From the data for $m=0$ (given in
\eqref{eq:m=0}) and for $m=1,2,3$, one is then easily
able to work out a (at this
point, conjectural) formula for the output of the algorithm, namely if
$$F(n,k)=\frac {(a+n)_k\,(b-n)_k\,(c)_k} {(a+n+1)_k\,(b+c+m-n)_k\,k!},$$
then
\begin{multline} \label{eq:rek2}
(b-n)_m\,(c)_m\,( a + n+1) (b-n-1)\\
\times (a - b - c-m + 2 n+1)
         (a - b - c-m + 2 n+2) F(n,k) \\+
(b-n)_m\,(c)_m\,      (a-b+2n+2) (a-b+2n+1)\kern5cm\\
\times (a-c+n+1)
       ( b + c+m - n-1) F(n+1,k] =\De_k F(n,k) R(n,k),
\end{multline}
where 
\begin{multline*}
R(n,k)=  \frac {k (b + c+m + k - n-1 )} {( a + n
       )  (  b + c +m- n-1 )  (  b + k - n-1 )}
\\\times
      ( k^{m+1} ( b + c +m- n-1 )  
          (  a + n +1)  \\
\cdot
          ( 1 - a - a^2 - 2 b + b^2 + a c + n - 2 a n - 2 b n +
            c n  -m(n+1-b) )\\  +\text{terms with lower powers in $k$}),
\end{multline*}
If one, for simplicity, replaces $b+c+m$ by $d$ in \eqref{eq:rek2},
sums both sides over $k$ from $0$ to $N$, and finally lets $N$ tend to
infinity, one arrives exactly at \eqref{eq:b11}, with $a$ replaced by
$a+n$, $b$ replaced by $b-n$, and $d$ replaced by $d-n$. 

As we pointed
out, this is at best a half rigorous derivation of \eqref{eq:b11} in the
case that the difference $d-b-c$ is a non-negative integer, but there
is no guarantee at all that this formula should also hold for {\it
any} $d$. (To explain two of the possible pitfalls: first, there 
are always two ways
to translate expressions such as $(c)_m$ into gamma functions:
$(c)_m=\Ga(c+m)/\Ga(c)=(-1)^m\Ga(1-c)/\Ga(1-c-m)$. These lead to
different formulae if $m$ is replaced by $d-b-c$, where $d$ is
arbitrary. Second, sometimes one may even miss whole additional terms
in a formula, which one does not see if some parameter is specialised
to a non-negative integer because this additional term happens to
vanish for this specialisation.) However, one can now {\it prove}
\eqref{eq:b11} continuing along the above lines: first, one verifies
that \eqref{eq:b11} is valid for $d=b$ by using Gau{\ss}' summation
formula \eqref{eq:Gausz}. Next, one replaces $d$ by $b+n$ in
\eqref{eq:b11}, and one uses the Gosper--Zeilberger algorithm to find
recurrences in $n$ for the left-hand and the right-hand sides of
\eqref{eq:b11}. Thus, one knows that \eqref{eq:b11} holds with $d=b+n$
for any non-negative integer $n$. Since both sides of \eqref{eq:b11}
are analytic in $d$ in a neighbourhood of $\infty$, one can use the
principle of analytic continuation to deduce that \eqref{eq:b11} holds
for any complex $d$ where both sides are defined.

We did that, but finally
we did succeed to work out a proof using known contiguous
relations. Since this is completely elementary and, as we believe, more
instructive, this is the proof that we have included 
in Section~\ref{sec:contig}. For
obtaining the general identities which are behind \eqref{eq:id2} and
\eqref{eq:id3}, given in
Propositions~\ref{prop:beta} and \ref{prop:gamma3}, we proceeded
similarly. In fact, the contiguous relations 
\eqref{eq:C55}, \eqref{eq:C15}, \eqref{eq:C27}, \eqref{eq:C54}, 
which we used in the proofs, are {\tt C55}, {\tt C15}, {\tt C27}, and
{\tt C54}, respectively, in HYP.

Our computer experiments suggest that the above procedure produces a
relation of the type \eqref{eq:b11} for {\it any} series
\begin{equation*}
{} _{3} F _{2} \!\left [ \begin{matrix} 
    a+a_1n,b+b_1n,c+c_1n\\ a+a_1n+a_2,d+d_1n
    \end{matrix} ;
    {\displaystyle 1}\right ],
\end{equation*}
as long as $a_1,a_2,b_1,c_1,d_1$ are integers, $a_2$ a positive
integer, and $b_1+c_1=d_1$. 
However, most of the time none of $a,b,c,d$ 
appears linearly in the big polynomial factor on the right-hand side.
This makes it difficult to extract a general solution of the
Diophantine equation which arises when one equates the polynomial
factor to zero. Nevertheless, experimentally, there are many solutions
for various choices of $a_1,b_1,c_1,d_1$.

\def\refname{Bibliography}

\end{document}